\documentclass[11pt,preprint]{article}
\date{August 1, 2016}
 \def\dateref{\hfil August 1, 2016 - preprint \hfil}
\usepackage [utf8]{inputenc}
\usepackage[math]{iwona}

\usepackage{amsmath,amsfonts, amssymb, version,  enumerate,color,draftcopy,hyperref}
\usepackage[amsmath]{ntheorem}
\usepackage{layout}
\definecolor{lightgrey}{rgb}{0.3,0.3,0.3}

\newtheorem{Theorem}{Theorem}[section]
\theoremstyle{change}
\newtheorem{Def}[Theorem]{Definition}
\newtheorem{Lemma}[Theorem]{Lemma}
\newtheorem{Proposition}[Theorem]{Proposition}
\newtheorem{Corollary}[Theorem]{Corollary}
\newtheorem{Rem}[Theorem]{Remark}

\newenvironment{Definition}{\begin{Def}\rm}{\end{Def}}

\newenvironment{Proof}{\noindent\textbf{Proof.}}{\null\hfill$\square$\medskip}
   \newenvironment{Remark}{\begin{Rem}\rm}{\end{Rem}}

\newcommand\half{\frac12}
\newcommand\irest[2]{\Bigl|_{#1}^{#2}\Bigr.}
\newcommand\norm[1]{\Vert #1 \Vert}
\newcommand\inv{^{-1}}
\newcommand\summ[3]{\displaystyle\sum_{#1 = #2}^{#3}}
\newcommand\BC[2]{\binom {#1}{#2}}

\newcommand\A{\mathbf{A}}
\newcommand\E{\mathbb{E}}

\usepackage{geometry}
\let\epsilon\varepsilon
\def\cD{{\cal D}}\def\cA{{\cal A}}\def\cB{{\cal B}}\def\cS{{\cal S}}
\let\ds\displaystyle

\def\colon{:\,}
\def\frak{\displaystyle\frac}
\def\tA{{\tt A}}\def\tB{{\tt B}}
\def\ta{{\tt a}}\def\tb{{\tt b}}
\def\N{{\bf N}}\def\Z{{\bf Z}}\def\R{{\bf R}}\def\w{{\bf w}}
\def\p{{\bf p}}\def\q{{\bf q}}
\def\e#1{{\bf e}_{#1}}\def\he#1{{\widehat{\bf e}}_{#1}}
\def\tA{{\widetilde A}}
\def\wtW{{\widetilde W}}
\let\wdtl\relax
\def\wtV{{\widetilde V}}
\def\crA{{\stackrel{\circ}{A}}}
\def\crB{{\stackrel{*}{B}}}
\def\crv{{\stackrel{\circ}{\bf v}}_\lambda}
\def\crw{{\stackrel{\circ}{\bf w}}_\lambda}
\def\hcrw{{\stackrel{\widehat\circ}{\bf w}}_\lambda}

\def\crz{{\stackrel{\circ}{\bf z}}_\lambda}
\def\ccz{{\stackrel{\circ}{\bf z}}}
\def\csz{{\stackrel{*}{\bf z}}}
\def\wtV{{\widetilde V}}

\def\v{{\bf v}}
\def\u{{\bf u}}
\def\z{{\bf z}}
\def\jj#1{{#1}^{(j)}}
\def\n{{\bf n}}\def\s{{\bf s}}\def\i{{\bf i}}\def\h{{\bf h}}
\def\ij{{(i,j)}}\def\j{{\bf j}}\def\k{{\bf k}}
\def\hw{{\widehat{\bf w}}}
\def\whi{{\widehat{\bf i}}}
\def\wih{{\widehat{i}}}
\def\whj{{\widehat{\bf j}}}
\def\wjh{{\widehat{j}}}
\let\wh\widehat
\def\hf{{\widehat{f}}}
\def\lrf#1{\|#1\|_r^\infty}
\def\brakt#1#2{\langle\,#1,#2\,\rangle}
\def\D#1#2{D_{#1|#2}}
\def\wM{\widetilde{M}}
\def\be{\begin{equation}}\def\ee{\end{equation}}
\def\bea{\begin{eqnarray}}\def\eea{\end{eqnarray}}
\def\prf{\medbreak\noindent{\bf Proof}:\enspace}
\def\qed{\hspace*{\fill}\hbox{\vrule height 7pt \kern-.3pt
    \vbox{\hrule width 7pt \kern6.6pt\hrule width 7pt }\kern-.3pt\vrule height 7pt
    }\par}
\let\vareps\varepsilon
\def\hateps{{\widehat\varepsilon}}
\def\li{\mathop{\rm li}_2\nolimits}

\def\sc#1{\em #1}
\numberwithin{equation}{section}
\makeatletter
  \def\footnotestar{\xdef\@thefnmark{}\@footnotetext}
\makeatother

\begin{document}

\title{Exhaustion of an interval by iterated Rényi parking}

%\address{School of Mathematics and Statistics, University College Dublin, Ireland}

\author{Michael Mackey and Wayne G. Sullivan}

%\email{wayne.sullivan@ucd.ie}
%\address{University College Dublin}

%\keywords{Rényi parking}

%\begin{keyword}
%   Rényi parking
% \end{keyword}

%\footnotestar{AMS Subject Classification }

\maketitle
\footnotestar{School of Mathematics and Statistics, University College
  Dublin\hfill mackey@maths.ucd.ie\\ \null\hfill wayne.sullivan@ucd.ie}

\begin{abstract}
  We study a variant of the Rényi parking problem in which car length
  is repeatedly halved and determine the rate at which the remaining
  space decays.
\end{abstract}

\pagestyle{myheadings}
\markboth{\dateref}{\dateref}

\section{Introduction}

Rényi modeled a long interval randomly jammed with cars of unit
length.   More precisely, each arriving car is parked at a location
chosen according to uniform probability distribution over the set of
available parking locations until no parking locations remain.  
The resulting gap lengths constitute a random variable with range $(0,1)$.  The constant
$C_R\approx 0.747598$  is the
limit, as the interval length tends to infinity, of the  proportion of
the interval covered and  was analytically
determined \cite{renyi}.  This Rényi parking constant is scale
invariant, in the sense that it holds for cars of any fixed length.
Variations on Rényi's parking problem consider higher dimensions, a
discretised setting, or cars of mixed lengths.  Modelling of physical
processes, such as random sequential adsorption, frequently refer to
the idealised Rényi model (see for example \cite{dana_cars, MR1824203,MR855759}
and references therein).  Here, we consider another variation. 
For our purposes it is more convenient to assume we have initially
jammed with cars of length $2$ leaving us with gaps distributed
according to a probability distribution $X_1$ with range $(0,2)$.  In
our formulation, the
next step is that these $X_1$-distributed intervals are jammed by cars
of length $1$ to give a new gap distribution $X_2$ with range
$(0,1)$.  We 
then introduce and jam with cars of length $\half$ to get $X_3$ with
range $(0,\half)$, then cars of length $\frac 14$ and so on.  It
is clear that our initial intervals are exhausted by this process.
\emph{Our motivating question is that of determining the rate of exhaustion.}
We will show that, if $L_n$ denotes the expected length of interval which
remains uncovered after $n$ stages, then $L_{n+1}/L_n$ tends to a
finite limit $R_\half\approx 0.61$.  We will numerically calculate this limit with
error bounds and also see that the limiting behaviour is largely
independent of the initial probability distribution $X_1$.  In
particular, $X_1$ need not be the Rényi distribution outlined above.
Proof of these facts requires finding a limiting distribution $X$ to
the sequence $\{X_n\}$ and for this we need to normalise at each stage.
Integral equations governing the evolution of $\{X_n\}$, their
corresponding densities $\{f_n\}$, and steady state analogues will be
found.   We approach the solution to these equations by considering
an infinite dimensional eigenvalue problem.  Numerical analysis with error
bounds shows that there is a unique eigenvalue of
maximum absolute value. Convergence of the iterates of the linear
system to the associated
eigenspace yields convergence of the sequence of probability densities.

The paper is laid out as follows.  In Section~2, we consider the sequence of cumulative
probability distributions  that arise as a result of the
physical process of jamming with cars of fixed length at each stage
before halving the car length for the next stage. This yields an
integral equation relating the distribution at Stage $n+1$ to that at Stage
$n$.   We deduce from this the evolution equation in the case of a
probability density
 which naturally leads to 
an integral equation governing any probability density  which is a fixed point of this process.
From this, we observe functional properties of any such fixed point 
which is sufficiently regular.

In Section~3, we use series representations of the iterated sequence of probability
densities in order to recast the integral formulation of density
evolution as an infinite linear system.  The sought-after asymptotic decay rate
$R_\half$  is seen to be one half of the spectral
radius of the infinite matrix  and
spectral properties of this matrix, $\A$, are explored.

Following mention of floating point concerns in Section~4, careful
numerical analysis of a similarity transformation of the matrix $\A$ is provided in Section~5, based on
Gershgorin's Theorem with error and truncation estimates.   Using the same
similarity transformation, we deduce 
a convergence result for iteration of $\ell^1$ vectors under $\A$,
corresponding to convergence of our probability densities.

Finally, in Section~6, we consider convergence for more general initial
distributions without density.

\section{The Model}
Let $N$ be the large number of gaps 
in a long interval of length $L$ formed, say, by parking cars of
length 2 until jamming occurs.   The gap lengths form a random
variable, $X_1$, with range $(0,2)$.   The expected number of these gaps which have
length greater than $1$ is
$NP(X_1>1)$.  Each of these has a car of length $1$ placed in it,
randomly, by uniform distribution over the set of available locations. 
This creates $2NP(X_1>1)$ new gaps, referred to as Stage 2 gaps, of length
in $(0,1)$ and removes $NP(X_1>1)$ Stage 1 gaps of length in $(1,2)$.  The
length of one of these $2NP(X_1>1)$ Stage 2 gaps is a random variable
with range $(0,1)$ which we call $Y$.  The distribution of all
$N+NP(X_1>1)$ gaps, whose length ranges over $(0,1)$, will be
represented by the random variable $\tilde X_2$ which we normalise to
$X_2=2\tilde X_2$ with range $(0,2)$.

We calculate the probability distribution of $\tilde X_2$.   Observe
that each gap at this point may have existed at the first stage or
have been created on introduction of cars at the
second stage.   We refer to these as Stage 1 gaps and Stage 2 gaps
respectively.  For $t\in
(0,1)$,  
\begin{equation}P(\tilde X_2<t) = P(\mbox{gap from Stage 1})P(X_1<t| X_1<1) + P(\mbox{gap from
  Stage 2})P(Y<t).\label{eq:der1}
\end{equation}
Clearly,
\begin{align*}
  P(\mbox{gap from Stage 1}) &= \frac{N-NP(X_1 >1)}{N+NP(X_1>1)} =
  \frac{1-C_1}{1+C_1} \\
\intertext{where $C_1=P(X_1>1)$ and}
P(\mbox{gap from Stage 2}) &= \frac{2C_1}{1+C_1}.
\end{align*}

We need $P(Y<t)$.  Consider  a Stage 1 gap of length $\lambda \in (1,2)$,
which has a unit length car randomly placed in it by way of uniform
distribution over $(0,\lambda-1)$, resulting in two new
Stage 2 (i.e. $Y$-)gaps.  The probability that one of these $Y$-gaps
has length $<t$,  $P(Y^{\lambda}<t)$ is $1$ if $\lambda<t+1$ and $\frac
t{\lambda-1}$ otherwise.  Hence
\begin{align*}
  P(Y<t) &= \int_{\lambda=1}^2 P'(X_1<\lambda | X_1>1)
  P(Y^\lambda < t) d\lambda \\
            &= \int_1^2 \frac{P'(X_1<\lambda)}{P(X_1>1)}
            P(Y^\lambda < t) d\lambda \\
            &= \frac 1{C_1} \int_{\lambda=1}^{t+1} 
            {P'(X_1<\lambda)} d\lambda + \frac 1{C_1}
            \int_{\lambda=t+1}^2 P'(X_1<\lambda) \frac
            t{\lambda-1} d\lambda \\
            &= \frac 1{C_1}\left[ P(X_1<t+1) -P(X_1<1) + \frac
              {tP(X_1<\lambda)}{\lambda-1}\irest{t+1}2 +
              t\int_{t+1}^2
              \frac{P(X_1<\lambda)}{(\lambda-1)^2} d\lambda
            \right] \\
            &= \frac 1{C_1}\left[ t-1+C_1+t\int_{t+1}^2 \frac
              {P(X_1<\lambda)}{(\lambda-1)^2} d\lambda\right] .\\
\end{align*}

%\textbf{Remark.}  Since $\lim_{t\to 0}P(Y<t)=0$, we have \[\lim_{t\to
%  0} t\int_{t+1}^2 \frac{P(X_1<\lambda)}{(\lambda-1)^2} d\lambda = {1-C}.\]

Substituting into \eqref{eq:der1} we have
\begin{align*}
  P(\tilde X_2 < t) &= \frac{1-C_1}{1+C_1}P(X_1<t | X_1<1) + \frac
  {2C_1}{1+C_1}.\frac 1{C_1}\left[t+C_1-1+t\int_{t+1}^2
    \frac{P(X_1<\lambda)}{(\lambda-1)^2} d\lambda \right]\\
  &=  \frac{1}{1+C_1}P(X_1<t) + \frac
  {2}{1+C_1}\left[t+C_1-1+t\int_{t+1}^2
    \frac{P(X_1<\lambda)}{(\lambda-1)^2} d\lambda \right]
\end{align*}
Now for $X_2=2\tilde X_2$ and $t\in (0,2)$ we have
\begin{align*}
  P(X_2<t)&= P(\tilde X_2 < \frac t2) \\
  &= \frac 1{1+C_1}P(X_1<\frac t2) + \frac 2{1+C_1}\left[ \frac t2+C_1-1 +
    \frac t2
    \int_{(t+2)/2}^2 \frac {P(X_1<\lambda)}{(\lambda-1)^2}
    d\lambda\right] \\
  &=\frac1{1+C_1}\left[P(X_1<\frac t2) + t+2(C_1-1) +
  t\int_{1+t/2}^2\frac{P(X_1 < \lambda)}{(\lambda-1)^2} d\lambda\right].
\end{align*}
Letting $F_s(t)=P(X_s<t)$ and $C_s=P(X_s>1)=1-F_s(1)$,
the map $S:F_1 \mapsto F_2$ above can be applied
repeatedly to gain a sequence of
probability distributions $\{F_s\}$ on $(0,2)$, with associated
constants $C_s=P(X_s>1)$.
%\footnote{TODO: show that $S$ is, in some
%sense, a contraction, and converges under iteration.}  

We seek a
limiting distribution which
is a fixed point of this process and are led to consider the evolution equation
\begin{equation}\label{eq:distev} (1+C_s)F_{s+1}(x) = F_s(x/2) +x +2(C_s-1) + x\int_{(x+2)/2}^2
\frac{F_s(y)}{(y-1)^2}dy. 
\end{equation}
Assuming differentiability of $F_1$, we have a corresponding evolution of probability
densities $\{f_s\}$,
\begin{equation}
  (1+C_s)f_{s+1}(x)= \half f_s(x/2)+ \int_{1+x/2}^2
  \frac{f_s(y)}{y-1}dy \label{eq:int_ev}
\end{equation}
and another differentiation yields the form
\begin{equation}
  \label{eq:diff_ev}
(1+C_s)x{f_{s+1}}'(x) = \frac x4 {f_s}'(\frac x2) - f_s(\frac{x+2}2)
\end{equation}
with $f_{s+1}$ fixed by the requirement that $\int_0^2 f_{s+1}=1$.  Note  that if $f_s$ is continuously differentiable, then so too is $f_{s+1}$.

%Equation \eqref{eq:diff_ss} has been reached more efficiently by WGS!
%[as an integral equation: renyip.tex (2)].   

Remark that if $F_s \to F$ pointwise, then the limit distribution $F$
will satisfy the steady state version of \eqref{eq:distev}, namely
\begin{equation}\label{eq:distss} (1+C)F(x) = F(x/2) +x +2(C-1) + x\int_{(x+2)/2}^2
\frac{F(y)}{(y-1)^2}dy
\end{equation}
where $C=1-F(1)$.    There are corresponding steady state equations for the limiting
probability density $f$ in both integral form, 
\begin{equation}
  (1+C)f(x)= \half f(x/2)+ \int_{1+x/2}^2
  \frac{f(y)}{y-1}dy \label{eq:int_ss}
\end{equation}
and, assuming it is differentiable, differential form,
  \begin{equation}
    \label{eq:diff_ss}
    4(1+C) xf'(x) = xf'(\frac x2) -4f(1+\frac x2), \qquad (x\in (0,2))
  \end{equation}
where $C=\int_{(1,2)} f$.  This is a linear first
order ordinary differential equation whose solution is hampered by the
presence of two time delays and non-linearity due to the
dependence of $C$ on the solution.    We will establish uniform convergence of
the sequence of probability densities on compact subsets of $(0,2]$, from a general starting
point, $f_0$, to a solution of \eqref{eq:int_ss} and pay  particular attention  to
numerical estimation of $C$.

\subsection{Rate of Decay of Remaining Space}
Consider the question of exhaustion rate of an interval by
this iterative process.  Let $X$ be the random variable representing
gap lengths whose distribution  $F$ 
satisfies \eqref{eq:distss} with $C=1-F(1)=P(X>1)$.   If $F$ is differentiable then we 
let $f=F'$ be the corresponding density function.  For a large
number $N$ of such $X$-gaps, the total length has
expected value $N\E(X)$.
Jamming with unit length cars, we can expect to fit $N P(X>1)$ of same
with total length $N P(X>1)$.  The uncovered length has reduced to
$N\E(X)-NP(X>1)$ and the ratio $R_\half$ is
$(N\E(X)-NP(X>1))/(N\E(X))$ or
\begin{equation}
  R_\half=1- \frac{P(X>1)}{\E(X)} = 1-\frac C{\E(X)}. \label{Rhalf}
\end{equation}
The expected value $\mathbb{E}(X)=\int_0^2 tF(dt)$ ($= \int_0^2
tf(t)\,dt$)  $= 2-\int_0^2 F$.  
On integrating \eqref{eq:distss} and changing order in the resultant
double integral we find
\begin{align*}
  (1+C)\int_0^2 F &= 2\int_0^1 F +2+4(C-1)+ \int_{x=0}^2\int_{y=1+\frac
    x2}^2 \frac {xF(y)}{(y-1)^2}\, dy\,dx \\
   &= 2\int_0^1 F + 4C-2+\int_{y=1}^2 2F(y)dy\\
   &= 2\int_0^2 F +4C-2.
\end{align*}
Thus $\int_0^2 F = \frac{2(1-2C)}{1-C}$ and $\mathbb{E}(X)=\frac{2C}{1-C}$.
Now, from \eqref{Rhalf}, we have the rate of exhaustion as
\begin{equation}
R_\half=\frac{1+C}2.
\label{decay}
\end{equation}

\subsection{Properties of a solution}

If $f$ is a continuous probability density on $(0,2)$ satisfying \eqref{eq:int_ss}
then $\lim_{x\to2} f(x)$ exists and we may take $f(2)$ to be this
limit.  Then, 
$f(1)=2(1+C)f(2)$ and, as $C$ is a probability, it follows \begin{equation}f(1)\le
4f(2)\le 2f(1).\label{eq:ineq}  \end{equation}
                                                           
\begin{Proposition}\label{prop:21}
  Let $f_0$ be a probability density on $(0,2]$, continuous in a
  neighbourhood of $0$, 
  generating the sequence of probability densities $\{f_s\}$ according
  to \eqref{eq:int_ev}.  There exists $\sigma\in \N$ such that
  $f_s>0$ on $(0,2]$ for all $s>\sigma$.  Indeed, for each $s>\sigma$,
  there exists $\epsilon_s >0$ such that $f_s> \epsilon_s$ on $(0,2]$.
\end{Proposition}

\begin{Proof}   By hypothesis, $f_0$ is continuous  on $(0,a)$ for
  some $a>0$ and so 
  \eqref{eq:int_ev} implies that $f_1$ is continuous on $(0,2a)$.  
After a finite number of iterations we
  have a continuous probability density on $(0,2)$.   Without loss of
  generality therefore, we may assume that $f_0$ is continuous, as are
  subsequent iterations.  
 Choose $\alpha\in(0,2)$
  such that $f_0(\alpha)>0$.  Indeed, by continuity, we may assume
  $\alpha$ is not a (negative) power of 2, and thus for some $k\in\N$,
  $\beta=2^k\alpha \in (1,2)$.  From \eqref{eq:int_ev},
  $f_1(2\alpha)\ge \frac14 f_0(\alpha)>0$, and by repetition,
  $f_k(\beta)>0$.   Continuity of $f_k$ implies that the integral
  $\int_{1+\frac x2}^2 \frac{f_k(y)}{y-1}dy >0$ whenever
  $1+\frac x2<\beta$, that is $x<2(\beta-1)$ and therefore
  $f_{k+1} > 0$ on $(0,2(\beta-1))$.  In particular, there exists
  $j\in\N$ such that $f_{k+1}>0$ on $(0,2^{-j})$ and it follows, as
  above, that for $\sigma=k+j+1$, $f_\sigma>0$ on $(0,2)$.    
  Remark also that
  $f_\sigma>0$ on $(0,2)$ implies  that $f_{\sigma+1}>0$ on $(0,2]$ .

Next, we can choose $\delta>0$ such that $f_\sigma(y)>\half f_\sigma(1)$ for
  $y\in(1,1+\delta)$.   Then $f_{\sigma+1}(x) \ge \frac 14 f_s(\frac
  x2)+\half\int_{1+\frac x2}^2
  \frac{f_\sigma(y)}{y-1}dy>\frac14f_\sigma(1)\int_{1+\frac x2}^{1+\delta}
  \frac{1}{y-1}dy \to\infty$ as $x\to0^+$.    That is, $\lim_{x\to
    0^+} f_{\sigma+1}(x)=\infty$.  In particular, there
  exists $x_0\in(0,2)$ such that $\inf_{(0,2]} f_{\sigma+1} =\inf_{[x_0,2]}
  f_{\sigma+1} > 0.$
\end{Proof}

\begin{Remark}
  The hypothesis in the above result can be weakened to deal with more
  general probability measures.  See Proposition~\ref{prop:mp}.
\end{Remark}

% \begin{Lemma}
%   Suppose $f_s>0$ on $(0,2]$.   Then there exists $\epsilon>0$ such
%   that $f_{s+1}\ge \epsilon$ on $(0,2]$.  
% \end{Lemma}

% \begin{Proof}  
%   Choose $\delta>0$ such that $f_s(x)>\half f_s(1)$ for
%   $x\in(1,1+\delta)$.   Then $f_{s+1}(x) \ge \frac 14 f_s(\frac
%   x2)+\half\int_{1+\frac x2}^2
%   \frac{f_s(y)}{y-1}dy>\frac14f_s(1)\int_{1+\frac x2}^{1+\delta}
%   \frac{1}{y-1}dy \to\infty$ as $x\to0^+$.    In particular, there
%   exists $x_0\in(0,2)$ such that $\inf_{(0,2)} f_{s+1} =\inf_{[x_0,2]}
%   f_{s+1} > 0.$
% \end{Proof}

\begin{Proposition}\label{prop:23}
  Suppose $f$ is a continuously differentiable probability density on $(0,2]$ which is a fixed point of \eqref{eq:int_ev} (that is, a solution of \eqref{eq:int_ss}).  
  \begin{enumerate}[(i)]
    \item $f$ is strictly positive.\label{prop:23i}
    \item If there
      exists $M\in\R$ such that $f'\le M$ on $(0,\epsilon)$ for some $\epsilon>0$ then $f$ is
      decreasing.\label{prop:23ii}
      \item If further, $f$ is twice continuously differentiable and there exists $N\in\R$ such that $f'' \ge N$
        then $f$ is convex.\label{prop:23iii}
  \end{enumerate}
\end{Proposition}

\begin{Proof}  The first part follows immediately from Proposition~\ref{prop:21}.
  For the second part, suppose, for sake of contradiction, that
    $f'(\xi)\ge 0$ for some $\xi\in(0,2)$.  Then, using
    \eqref{eq:diff_ss},
    \[0 \le \xi f'(\frac \xi2) - 4f(1+\frac \xi2)\]
    and so $0<2f(1+\frac\xi2 ) \le \frac\xi2 f'(\frac\xi2)$.  In
    particular, $f'(\frac\xi2) > 0$ and we can repeat the argument
    to gain, for every $n\in\N$,
\[ 0< 2f(1+\frac\xi{2^n}) \le \frac\xi{2^n} f'(\frac\xi{2^n}).\]
Thus $\limsup_n \frac\xi{2^n}f'(\frac\xi{2^n}) \ge 2f(1)>0$ which is
incompatible with $f'$ being bounded above near $0$, proving the
claim. 

Towards the statement on convexity, we may differentiate \eqref{eq:diff_ss} to find for all $x\in (0,2)$
\begin{align*} K(f'(x)+xf''(x))&= f'(\frac x2) +\frac x2f''(\frac x2)-2f'(1+\frac x2) \\
  &= f'(\frac x2) -4\frac {f(1+\frac x2)}x + 4\frac {f(1+\frac x2)}x + \frac x2f''(\frac x2) -2f'(1+\frac x2)\\
  Kxf''(x)&= 4\frac {f(1+\frac x2)}x + \frac x2f''(\frac x2) -2f'(1+\frac x2), 
\end{align*}
where $K=4(1+C)$.  Suppose $f''(\xi)<0$ some $\xi\in (0,2)$.  Then the above, together
with positivity of $f$ and $-f'$, implies $\frac \xi2 f''(\frac
\xi2)<K\xi f''(\xi)<0.$  We can repeat with $\frac\xi2$ in place of $\xi$ to get \[\frac \xi4 f''(\frac \xi4) <   K\frac\xi2 f''(\frac\xi2) < K^2 \xi f''(\xi).\]
Further repetition yields, $f''(\frac\xi{2^n}) < (2K)^n f''(\xi) \to -\infty$ which does not allow that $f''$ is bounded below.
\end{Proof}

If $\{f_n\}$ is a sequence of $\mathcal{C}^1$ functions which converge
to $f$, also $\mathcal{C}^1$, then it does not follow that $f_n'\to
f'$, even if the convergence is uniform--  a
common example provided is the sequence $f_n(x)=\frac 1n \sin
(nx) \to 0$ on $[0,\pi]$, which even has a uniform bound on
$\{\norm{f_n'}\}$.   Nevertheless, an application of the mean
value theorem yields the following.

\begin{Lemma}\label{eq:bdd}
    Suppose $\{f_n\}$ is a sequence of $\mathcal{C}^1$ functions which converge
to $f$, also $\mathcal{C}^1$, on an interval $I$.  If $\sup_{n,x\in I}
|f_n'(x)| \le M$ then $\sup_{x\in I} |f'(x)| \le  M$.
\end{Lemma}

\newcommand\eps{\epsilon}
% \begin{Proof}
% Suppose, for sake of contradiction, that $\sup_I|f'(x)|\ge M+\eps$ for
% some $\eps >0$.    Then we may choose  a subinterval $[\alpha, \beta]$
% of $I$ such that $|f'(x)|\ge M+\frac\eps2$ for all
% $x\in[\alpha,\beta]$.  The Mean Value Theorem implies
% $\absval{\frac{f(\beta)-f(\alpha)}{\beta-\alpha}}\ge M+\frac\eps2$ .  Then by
% convergence of $(f_n)_n$,  choose $n_0$ such that
% $\absval{f_{n_0}(\beta)-f_{n_0}(\alpha)} \ge
%   (M+\frac\eps4)(\beta-\alpha)$.   Again, the Mean Value Theorem
%   implies that $|f_{n_0}'|\ge M+\eps/4$ at some point which is a contradiction.
% \end{Proof}

\begin{Corollary}\label{cor:bdd}
    Suppose $\{f_n\}$ is a sequence of $\mathcal{C}^1$ functions which converge
to $f$, also $\mathcal{C}^1$, on $I=(0,L]$.  If $\sup_{n,x\in I}
|x f_n'(x)| \le M $ then $\sup_{x\in I} |xf'(x)| \le  M$.
  \end{Corollary}

  \begin{Proof}
    Let $g_n(x)=xf_n(x)+\int_x^L f_n(y)dy$.   For $x\in (0,L)$, we
    have uniform convergence of $f_n$ to $f$ on the compact interval
    $[x,L]$ and it follows $\int_x^L f_n(y)\,dy \to \int_x^L
    f(y)\,dy$.  Then $g_n(x) \to g(x):=xf(x)+\int_x^L f(y)\,dy$ for all
    $x\in (0,L]$.  Since $g_n'(x)=xf_n'(x)$ is bounded by hypothesis,
    and $g'(x)=xf'(x)$,
    Lemma~\ref{eq:bdd} guarantees that $\sup_{x\in(0,L]} |xf'(x)|$ has
    the same bound.
  \end{Proof}

\begin{Proposition}\label{prop:props}
  Let $f_0$ be a continuously differentiable probability density on
  $(0,2]$ with $\sup_{(0,2)} |xf_0'(x)|<\infty$,  and let
  the sequence $\{f_s\}$ be generated by \eqref{eq:int_ev}.
%\[ (1+C_s)f_{s+1}(x) = \half
%  f_s(\frac x2) + \int_{1+\frac x2}^2 \frac{f_s(y)}{y-1}\, dy\]
%  where $C_s= \int_1^2 f_s$.  
If $f_s\to f\in \mathcal{C}^1((0,2])$  pointwise  then
  \begin{enumerate}[(i)]
    \item $f$ is a probability density and a solution of \eqref{eq:int_ss},
      \item $f>0$,
        \item $f$ is decreasing,
  \item $\sup_{(0,2]} |tf'(t)|$ is finite,
    \item $\lim_{t\to0} tf'(t)$ exists in $(-2,0)$.
  \end{enumerate}
\end{Proposition}

\begin{Proof}
%\item[(i)-(iii)] 
\begin{enumerate}  
\item[(iv)]
  Each $f_s$ is continuously differentiable on $(0,2]$,
  and \begin{equation}K_s xf_{s+1}'(x)= x f_s'(\frac
  x2)-4f_s(1+\frac x2)
\label{eq:qq}\end{equation}
where $K_s=4(1+C_s) \ge 4$.  

 Since $f_s\to f$ uniformly on the compact set $[1,2]$,
 $D:=\sup_{\{x\in[1,2],s\in \N\}} f_s(x) <\infty$.
Letting $g_s(x)=xf'_s(x)$, equation \eqref{eq:qq} reads \[g_{s+1}(x) =
\frac2{K_s}g_s(\frac x2)  - \frac 4{K_s}f_s(1+\frac x2)\]
and \[ |g_{s+1}(x)|  \le 
\half |g_s(\frac x2)|+D.\]
Let $G_s=\sup_{(0,2]}
|g_s(x)|$.  Then $G_{s+1} \le \half G_s +D$ for each $s$.  By
hypothesis, $G_0$ is finite. and it follows that every $G_s$ is finite
with \[G_s \le \frac 1{2^s}G_0 + (1+\half+\cdots+\frac 1{2^{s-1}})D \le
G_0+2D.\] 
  That is, there exists
$M>0$, such that for all $s\in
\N$, $\sup_{(0,2]} |xf_s'(x)| < M$ and it follows from \ref{cor:bdd} that $\sup_{(0,2]} |xf'(x)|\le
M<\infty$.

\item[(i)]  From the estimates above, $\int_a^b |f_s'(x)| dx \le
  |\int_a^b \frac Mx \, dx| \le M|\log(b/a)|$.  Then $f_s(a)\le
  f_s(b)+M|\log(b/a)|$.   For each $s$, we may choose $b\in [1,2]$ with
  $f_s(b)\le 1$ so that $f_s(a)\le 1 + M\log(2/a)$ for all
  $a\in(0,1]$.  Then for $0<\delta\le1$,
\[\int_0^\delta f_s(x)\,dx \le \delta(1 + M(1+\log 2 ))
  +M\delta\log \delta\inv.  \]
Hence, as $\delta\to 0$, $\int_0^\delta f_s \to 0$, and
$\int_{\delta}^2 f_s(x)\,dx \to 1$ uniformly in $s$.

Since $f_s\to f$ uniformly
  on $[1,2]$, we have convergence of the integrals $C_s=\int_1^2 f_s$
  to $C=\int_1^2 f$, and of $\int_{1+x/2}^2 \frac{f_s(y)}{y-1}dy$ to
  $\int_{1+x/2}^2 \frac{f(y)}{y-1}dy$ for every $x>0$.  Thus we take
  limits in \eqref{eq:int_ev} to conclude $f$ satisfies
  \eqref{eq:int_ss}.  It is clear that $f$ is non-negative valued and
  on considering compact sub-intervals $[\delta,2]$, we also
  have that $\int_{0^+}^2f=\lim_{\delta \downarrow 0} \int_\delta^2  f = \lim_{s} \int_{0^+}^2 f_s = 1$.

\item[(ii)]  follows from (i) and  Proposition~\ref{prop:23}.

\item[(v) and (iii)] 
Letting $L=\limsup_{x\to 0^+} xf'(x)$, equation \eqref{eq:diff_ss}
  implies $4(1+C)L=2L-4f(1)$. From (iv), we know $L$ is finite and
  from (ii) $f(1)>0$ so
  $L(4C+2)=-4f(1)<0$ and $L=\frac{-4f(1)}{4C+2}<0$.  Notice
  the same quantity arises on taking the limit infimum (which is also
  finite from (iv)), and
  thus $\lim_{x\to 0^+} xf'(x)=L=\frac{-4f(1)}{4C+2}\in (-\infty,0)$.
  This implies that $f'$ is negative in a neighbourhood of $0$ and so from
  Proposition~\ref{prop:23}(\ref{prop:23ii}), $f$ is decreasing on $(0,2]$.
  Moreover, since $f$ is decreasing and \eqref{eq:ineq} $f(2)\ge f(1)/4$,
  $1=\int_0^2 f \ge f(1)+f(2) \ge \frac 54 f(1)$.  Hence $f(1)\le 4/5$
  and $L=-\frac {4f(1)}{4C+2} \ge -\frac 85$.

% \item[(ii)] Since $f$ is decreasing, it suffices to show that
%   $f(2)>0$.  The alternative is that $f(2)=0$ and hence by
%   \eqref{eq:int_ss} $f(1)=0$ and by monotonicity, $f|_{[1,2]}\equiv
%   0$ so $C=0$.   In this case, the integral partof \eqref{eq:int_ss}
%   contributes nothing and we see $f(2^{-n})=0$ for all $n$.  Again
%   monotonicity implies that $f|_{(0,2]}\equiv 0$ which contradicts the
%   fact that $f$ is a probability density.   
\end{enumerate}
\end{Proof}

Later (Corollary~\ref{fstarcor}), we prove the existence of a probability density
$f^*$ on $(0,2]$ of the form $f^*(t)=\ell^*\log (t)+\sum_{k=0}^\infty
a_k^*(t-1)^k$ which is a solution of \eqref{eq:int_ss}.   The
coefficient vector $(\ell^*,a_0^*,a_1^*,\ldots)$ is in $\ell^1_\rho$ for
$1\le \rho\le 3$ (see Definition~\ref{def:rho}) .  This means $f^*$ is a $\mathcal{C}^\infty$ function
on $(0,2]$.  Moreover, $f^*$  is
reached by iterating \eqref{eq:int_ev} from any initial probability
density $f_0=\ell \log(t) + \sum_{k=0}^\infty a_k(t-1)^k$ with $(\ell,
(a_k)_k) \in
\ell_1$, for example, from $f_0(t)\equiv \half$.  Consequently,
$f^*$ satisfies the conclusions of Proposition~\ref{prop:props}.  In
particular, we have from Proposition~\ref{prop:props}(v)that
$\ell^*=\lim_{t\to 0} t{f^*}'(t) < 0$ which, in turn, guarantees that ${f^*}''$ has a finite infimum over
$(0,2]$, so by Proposition~\ref{prop:23}(\ref{prop:23iii}), $f^*$ is also convex.

\iffalse
\begin{Corollary}
Let $f(t)$ be a probability density on $(0,2]$ of the form
\be
   f(t)=\ell\log(t)+\sum_{k=0}^\infty a_k(t-1)^k
\ee
with iterates given by $f_s(t)$ by (\ref{eq:int_ev}). Then
for each $0<a<2$ 
\be
   \lim_{s\to\infty} \sup_{a\leq t\leq 2}|f_s(t)-f^*(t)|=0,
\ee
where $f^*(t)$ is the probability density of the form
\be\label{fstardef}
   f^*(t)= c\Big(\v_0+\sum_{i=1}^\infty \v_i(t-1)^{i-1}),
\ee
where $\v$ is the eigenvector of $\A$ corresponding to the largest
eigenvalue  and the constant $c$ is chosen to so that
$\int_0^2 f^*\,dx=1$.
make $f^*$ a probability density.

\end{Corollary}
\fi

% \textbf{Remark.} The above result belongs earlier in the section about
% properties of $f_s$ and $f$ which satisfy (\ref{eq:int_ss}). One can
% refer to the Theorem above to be proved later.

\section{Linearisation}
We write successive iterations of an initial probability density $f_0$
on $(0,2]$
as
\[f_s(t)= \ell_s\log t + \summ k0\infty a_{k,s} (t-1)^k\]
with the desired limit function written as \begin{equation}f(t)
  =\ell\log t + \summ k0\infty a_k(t-1)^k.\label{eq:fform}\end{equation} 
% The use of the log term in the series expansion is motivated by the finite
% limit of $tf'(t)$ as $t$ approaches $0$, a fact established in
% Proposition~\ref{prop:props}.
We will identify $f_s$ with the sequence of coefficients $(\ell_s, a_{0,s},
a_{1,s},\ldots)$.  

Substituting the series representation of $f$ into \eqref{eq:int_ss} we have
\begin{align*}\int \frac{f(\tau)}{\tau-1}d\tau &= \ell\int \frac{\log
  (1+(\tau-1))}{\tau-1}d\tau + \int \sum_{k\ge 1} a_k(\tau-1)^{k-1}d\tau +
\int \frac{a_0}{\tau-1}d\tau \\
  &= \ell \int \sum_{k\ge 1} \frac{(-1)^{k+1}}k(\tau-1)^{k-1}\,d\tau +
  a_0 \log(\tau-1) + \sum_{k\ge 1}\frac{a_k}k(\tau-1)^k
\end{align*}
and
\begin{align*}
  \int_{1+\frac t2}^2  \frac{f(\tau)}{\tau-1}d\tau &=
  \ell\left(\sum_{k\ge 1} \frac{(-1)^{k+1}}{k^2}
    -\sum_{k\ge1}\frac{(-1)^{k+1}t^k}{k^22^k}\right)\\ &\qquad -a_0\log\frac t2  + \sum_{k\ge 1} \frac{a_k}k - \sum_{k\ge 1} \frac{a_kt^k}{k2^k}.
\end{align*}
Since $\sum_{k\ge1}\frac{(-1)^{k+1}}{k^2} =\frac{\pi^2}{12}$ and
$t^k=(1+(t-1))^k=\summ j0k \BC kj(t-1)^j$ we gain
\begin{align*}
  \int_{1+\frac t2}^2  \frac{f(\tau)}{\tau-1}d\tau &= \ell\left(
    \frac{\pi^2}{12} + \sum_{k\ge 1}\frac{(-1)^k}{2^k k^2}\summ j0k\BC
    kj(t-1)^j\right) - a_0\log \frac t2 + \sum_{k\ge 1} \frac{a_k}k
  -\sum_{k\ge 1} \frac{a_k}{k2^k}\summ j0k \BC kj(t-1)^j.
\end{align*}
We also find
\[ \half f(\frac t2) = \frac \ell2(\log \frac t2) + \sum_{k\ge0}
\frac{a_k}{2^{k+1}}\summ j0k\BC kj(t-1)^j(-1)^{k-j}\]
and thus \eqref{eq:int_ss} becomes
\begin{align*}
  (1+C)\left(\ell \log t + \sum_{k\ge 0}a_k(t-1)^k\right) & =(\frac \ell2-a_0)(\log \frac t2) + \sum_{k\ge0}
\frac{a_k}{2^{k+1}}\summ j0k\BC kj(t-1)^j(-1)^{k-j} \\
&\qquad +\ell\left(
    \frac{\pi^2}{12} + \sum_{k\ge 1}\frac{(-1)^k}{2^k k^2}\summ j0k\BC
    kj(t-1)^j\right) \\ &\qquad\quad+ \sum_{k\ge 1} \frac{a_k}k
  -\sum_{k\ge 1} \frac{a_k}{k2^k}\summ j0k \BC kj(t-1)^j.
\end{align*}
Equating the coefficients of $\log t$ gives $(1+C)\ell = \frac \ell2
-a_0$ while the constant terms yield \[(1+C)a_0 =(\frac\ell2
-a_0)(-\log 2) +\sum_{k\ge0} (-1)^k\frac{a_k}{2^{k+1}} +
\ell\left(\frac{\pi^2}{12}+\sum_{k\ge 1}\frac{(-1)^k}{k^22^k} + \sum_{k\ge
  1}\frac{a_k}k(1-2^{-k})\right) .\]  
The sum $\sum_{k\ge1} \frac{(-1)^k}{k^22^k}$ equals
$\operatorname{Li_2}(\frac32)$, approximately $-0.4484$,
where  $\operatorname{Li}_2$ denotes the
dilogarithm.  Finally, the coefficient of $(t-1)^r$, $r\ge 1$, gives
\[(1+C)a_r = \sum_{k\ge r} \frac{a_k}{2^{k+1}}\BC kr(-1)^{k-r} +
\ell\sum_{k\ge r} \BC kr\frac{(-1)^k}{2^kk^2} - \sum_{k\ge r} \BC kr
\frac{a_k}{k2^k}.\]
  The sum $\summ kr\infty \frac
{(-1)^k\BC kr}{k^22^k}$ is numerically slow to estimate but can be
simplified.

\begin{Lemma}
  Let $x\in[0,\half)$.  Then \[  
  \summ kr\infty \frac{(-1)^k\BC kr}{k^2}\left(\frac
    x{1-x}\right)^k = \frac{(-1)^r}{r}\summ kr\infty \frac 1k x^k .\]
In particular, taking $x=\frac 13$, \[ \summ kr\infty \frac
{(-1)^k\BC kr}{k^22^k} = \frac{(-1)^r}r \summ kr\infty \frac 1{k 3^k}.\]
\end{Lemma}

\begin{Proof}
  First, we use the expansion \(\log(1-x)= -\summ k1\infty
  \frac {1}k x^k\) to write
\begin{align*} \summ k1\infty \frac{1}k x^k &= \log
  \frac1{1-x}=\log (1-\frac {-x}{1-x}) 
  = -\summ k1\infty  \frac{(-1)^k \BC k1}{k^2}
\left(\frac x{1-x}\right)^k.
\end{align*}

This gives the desired identity when $r=1$ and we proceed by
induction; assuming the statement is true for $r$ we have,
\[ \summ kr\infty \frac {x^k}k =
  (-1)^r r\summ kr\infty \frac{(-1)^k\BC kr}{k^2}\left(\frac
    x{1-x}\right)^k.\]
Divide across by $x^r$,
% to gain
%\[\summ kr\infty \frac 1k x^{k-r} =
%  (-1)^r r\summ kr\infty \frac{(-1)^k\BC kr}{k^2}\frac
 %   {x^{k-r}}{(1-x)^k}\]
%and then 
differentiate, 
%\[\summ k{r+1}\infty \frac {(k-r)}k x^{k-r-1} =
%  (-1)^r r\summ kr\infty \frac{(-1)^k\BC kr}{k^2}\frac
 %   {(1-x)^k(k-r)x^{k-r-1}+x^{k-r}k(1-x)^{k-1}}{(1-x)^{2k}}\]
and multiply by $x^{r+1}$ to gain
\[\summ k{r+1}\infty x^k - r\summ k{r+1}\infty \frac {x^{k}}k =
  (-1)^r r\summ kr\infty \frac{(-1)^k\BC kr}{k^2}\left( \frac
    {kx^{k+1}}{(1-x)^{k+1}} + \frac{(k-r)x^k}{(1-x)^k}\right)
\]
which gives
\begin{align*}
\frac{x^{r+1}}{1-x}- r\summ k{r+1}\infty \frac {x^{k}}k 
%&=
%  (-1)^r r\summ kr\infty \frac{(-1)^k\BC kr}{k}\left[(\frac
 %   {x}{1-x}\right)^{k+1} \\ 
 %  &\qquad\qquad + (-1)^rr\summ k{r+1}\infty \frac{(-1)^k\BC
 %   kr(k-r)}{k^2}\left[  \frac{x}{1-x}\right]^k\\
&=
  (-1)^r \summ {k-1}{r-1}\infty (-1)^k\BC {k-1}{r-1}\left( \frac
    {x}{1-x}\right)^{k+1} \\ 
   &\qquad\qquad + (-1)^rr\summ k{r+1}\infty \frac{(-1)^k\BC
    kr(k-r)}{k^2}\left(\frac{x}{1-x}\right)^k.
\end{align*}
The first term on the right hand side can be simplified using \(\frac
1{(1-t)^{s+1}}= \sum_{k={s}}^\infty  \BC k{s} t^{k-s}\) leading to
\begin{align*}
\frac{x^{r+1}}{1-x}- r\summ k{r+1}\infty \frac {x^{k}}k 
%&=
%   (-1)^{r+1} \summ {k-1}{r-1}\infty (-1)^{k+1}\BC {k-1}{r-1}\left[ \frac
%     {x}{1-x}\right]^{k-1-(r-1)} \left(\frac x{1-x}\right)^{r+1}\\ 
%    &\qquad\qquad + (-1)^rr\summ k{r+1}\infty \frac{(-1)^k\BC
%     k{r+1}(r+1)}{k^2}\left[  \frac{x}{1-x}\right]^k\\
% &=
%   \left(  \frac x{1-x} \right)^{r+1} 
%    \summ {k-1}{r-1}\infty 
%        \BC{k-1}{r-1}  \left( \frac{-x}{1-x}\right)^{k-1-(r-1)}\\
%    &\qquad\qquad + (-1)^r  r\summ k{r+1}\infty  \frac{(-1)^k \BC
%     k{r+1}(r+1)}{k^2}\left[  \frac{x}{1-x}\right]^k\\
&=  \left(\frac x{1-x}\right)^{r+1} \hskip-1em \frac 1{ (1+\frac x{1-x})^{r}}
%\\ &\qquad\qquad 
+ (-1)^r  r(r+1)\summ k{r+1}\infty  \frac{(-1)^k \BC
    k{r+1}}{k^2}\left( \frac{x}{1-x}\right)^k.\\
\end{align*}
 Removing the common term and dividing by $-r$ gives the statement for
 $r+1$ as required.
\end{Proof}

The calculations above govern the relationship between the
coefficients of steady state solution $f$, but the corresponding
equations relating the coefficients of $f_{s+1}$ to those of $f_s$
amount to only a notational change and are summarised in the 
the following infinite system.
\begin{equation}
  \label{eq:lin}
  (1+C_s) \begin{bmatrix}\ell_{s+1}\\ a_{0,s+1}\\ a_{1,s+1}\\ \vdots
  \end{bmatrix}
=
\begin{bmatrix}
  \half &-1 &| &0 &0 &0 &\cdots \\
  0.027 &1.193 &| &\cdots &\frac 1k +\frac{(-1)^k}{2^{k+1}}-\frac
  1{k2^{k}} \\
-- &-- &| &-- &-- &-- &-- \\
\vdots & 0 & | & & & & \\
\frac {(-1)^r}r \summ jr\infty \frac 1{j3^j} &0 &| & & \ds\BC kr\left(
  \frac{(-1)^{k-r}}{2^{k+1}} -\frac 1{k2^{k}}\right) & & \\
\vdots &\vdots & | 
\end{bmatrix}
\begin{bmatrix}\ell_s\\ a_{0,s}\\ a_{1,s}\\ \vdots
  \end{bmatrix}
\end{equation}
where $k,r\ge 1$ and $\BC kr=0$ for $k<r$.  The constants 0.027 and 1.193 are just numerical placeholders for the
true values $\frac{\pi^2}{12}+\operatorname{Li_2}(\frac32)-\frac{\log 2}2$
and $\half+\log(2)$ respectively.   We have used $r$ and $k$ above to
index rows and columns respectively, but these indices are given relative to
the sub-matrices depicted.   This matrix of coefficients will hereafter be denoted by
$\A$ and its $(n+1)\times (n+1)$ principal submatrix by $\A_n$.  For example,
$\A_4$ is given by 
\[ \left[ \begin {array}{ccccc}  0.5&- 1.0& 0.0& 0.0& 0.0
\\ \noalign{\medskip} 0.02747926& 1.1931472& 0.25000000& 0.50000000&
 0.22916667\\ \noalign{\medskip}- 0.40546510& 0.0&- 0.25000000&-
 0.50000000& 0.062500000\\ \noalign{\medskip} 0.036065890& 0.0& 0.0&
 0.0&- 0.31250000\\ \noalign{\medskip}- 0.0055254100& 0.0& 0.0& 0.0&
 0.020833333\end {array} \right] 
\]

%\marginpar{Will we present this instead without the normalisation
%  constant, to give a genuinely linear system?}

Where necessary, $\A_n$ may be padded by zeroes so that expressions
such as $\A-\A_n$ and $\A_n-\A_m$ are well-defined.

\subsection{Spectrum}
Let us make some observations on $\A$.  As an operator on
$\ell^\infty$, $\A$ is not bounded since the absolute row sum of the
second row majorises the harmonic series.  $\A$ is however bounded on
$\ell_1$, with $\norm{\A}_1\simeq2.193$, as provided by the absolute second
column sum.  Indeed, for $k\ge 1$, the absolute sum of column $k$ is
convergent to $\half$ and bounded above by 1, while the absolute sum
of the only infinite column is bounded by 1.08, as provided by the
first of the following estimates.

\begin{Lemma}\label{li2bounds}
  \begin{align}
  (i)& \qquad
     \summ r1\infty \frac1r \summ jr\infty \frac1{j3^j} < 0.55, \\
    (ii)& \qquad  \summ rm\infty \frac1r \summ jr\infty \frac1{j3^j} <
      \frac 9{4(m^2 3^m)}.\label{col0bnd}
  \end{align}
\end{Lemma}

\begin{Proof}
  As $\summ jr\infty \frac1{j3^j}<\frac 1r\summ jr\infty
  \frac1{3^j}=\frac 32\frac1{r3^r}$,  the first double sum is
  bounded by $\frac32 \summ r1\infty
  \frac1{r^2 3^r}=\frac32\operatorname{Li_2}(2/3)\approx
  0.5493$.
The second double sum is similarly bounded by $\frac32 \summ rm\infty
\frac 1{m^2}\frac 1{3^r} = \frac 94\frac1{m^2 3^m}$.
\end{Proof}

While $\A$ is bounded on $\ell_1$, it is apparent that $\norm{\A_n
  -\A}_1 \to \half$, that is, $\A$ is not approximated by its finite
rank truncations.  For this reason, we will consider $\A$ acting on
the scaled $\ell^1$ space, $\ell^1_\rho$.  

\begin{Definition}\label{def:rho}
  For $\rho\ge 1$, $\ell^1_\rho=\{(x_n)_n: \norm{(x_n)_n}_{1,\rho}=\sum \rho^n|x_n| <\infty\}$.
\end{Definition}

If $D_\rho$ is the diagonal matrix $\operatorname{diag}(1,\rho,\rho^2, \ldots)$ then
we will consider $\A^\rho:=D_\rho \A D_\rho\inv$, that is, $(\A^\rho)_{i,j}=
(\A)_{i,j}.\rho^{i-j}$.   The norm of $\A$
as a bounded linear operator on $\ell^1_\rho$ is given by 
$\norm{\A}_{\ell^1_\rho}= \norm{\A^\rho}_{\ell^1}$.  In particular,
$\A \in \mathcal{L}(\ell^1_\rho)$ for $1\le\rho \le 3$
since the absolute column sums of $\A^\rho$ converge to zero and the
first column sum is finite for $\rho \le 3$.  Moreover,
$\norm{\A^\rho - \A^\rho_n}_1=\norm{\A - \A_n}_{\ell^1_\rho}\to 0$
which implies that $\A$ is approximated by its finite rank truncations
on $\ell^1_\rho$ for $1<\rho\le 3$ and so is a compact operator there.
Thus, its
spectrum on  $\ell^1_\rho$, $1<\rho\le 3$ consists only of
non-zero eigenvalues with a possible cluster point at zero. 
Moreover,  we have convergence of the spectrum of $\A_n$ to that of
$\A$ (see for example \cite[XI-9.5]{DS2}).   Of particular
interest is the eigenvalue corresponding to the spectral radius (we
will see in due course that there is only one such, and it is positive)
and the corresponding left and right eigenvectors.  Our immediate
objective is to estimate these to within acceptable tolerances.

\subsection{Similarity Transforms}

The first step is to find a similarity transform whose action on $\A$
yields a matrix for which
the application of Gershgorin's Theorem gives a good bound for the
maximum eigenvalue. This is done by choosing  approximate row and
column eigenvectors for $\A_m$. We express these here as $\hw$ and $\v$.
Let the standard unit column and row
vectors be denoted, respectively,  by $\e i$ and $\he i$, $\,i=0,\ldots,m$.
For a similarity transform we require that $\hw \v =1$ and the $j^{\rm th}$
entry of $\v$ satisfies $\v_j=\he j \v= 1$.
 The identity operator
$I$ and the operator $I^{(j)}$ which has the projection $\e j \he j$
removed are given by
\be
  I:=\sum_{i=0}^m \e i \he i,\quad 
  I^{(j)}:=\sum_{i\neq j}^m \e i \he i.
\ee
We use the similar notation $\v^{(j)}$ and $\hw^{(j)}$ to denote the
vectors with the $j^{\rm th}$ entry set to 0:
$\hw\jj I=\jj \hw$, $\jj I \v=\jj \v$.
For fixed $j$, $0\leq j\leq m$, we define two matrices
\be\label{wvdef}
  \wdtl W:=\e j\hw-\v\he j+I=\e j\hw-\jj\v\he j+\jj I,
  \quad\wdtl V:=\v(\he j-\jj\hw )+\jj I.
\ee
\begin{Lemma}
  For fixed $j$, let $\he j \v=1$ and $\hw\,\v=1$. Then
  $\wdtl W,\,\wdtl V$ as given by (\ref{wvdef}) 
  satisfy $\wdtl V\wdtl W=\wdtl W\wdtl V=I$. If, in addition,
  $\A\,\v=\lambda\v$ and $\hw\, \A=\lambda\hw$,
  then the entries of $\wdtl W \A\,\wdtl V$ are given by
\be\label{lessmbnd}
\begin{cases} \big( \wdtl W\A\,\wdtl V\big)_{jj}=\lambda &{} \\
       \big(\wdtl W\A\,\wdtl V\big)_{ij}=
      \big(\wdtl W\A\,\wdtl V\big)_{ji}=0& \text{if $i\neq j$};\\
 \big( \wdtl W\A\,\wdtl V\big)_{ik}=(\A)_{ik}-\v_i(\A)_{jk}
   & \text{otherwise.} 
\end{cases}
\ee
\end{Lemma}

\begin{Proof} $(\e j \hw)\wdtl V=\e j(\he j-\jj\hw)+\e j\jj\hw$.
$(- \v\he j)\wdtl V=(-\v_j)\v (\he j-\jj\hw )+0$. Also
$I\,\wdtl V=\v(\he j-\jj\hw ) +\jj I$. Collecting terms using $\v_j=1$
we find $\wdtl W\wdtl V=\e j\he j+\jj I=I$. 
A similar argument shows $\wdtl V\wdtl W=I$.
For the second part, $\A\wdtl V=\lambda\v(\he j-\jj\hw )+\A\jj I$.
Then using $\hw \v=\he j\v=1$,
\be\label{eq:46}
  (\e j \hw - \v\he j+I )\lambda \v(\he j-\jj\hw )=
  \lambda (\e j -\v +\v) (\he j-\jj\hw )=
  \lambda  (\e j\he j-\e j\jj\hw)  .
\ee
\be
  (\e j \hw -\jj \v\he j+\jj I )\A\jj I=
  \lambda \e j \jj\hw -\jj\v (\he j \A\jj I)+\jj I\A\jj I.
\ee
\be\label{eq:47}
  \wdtl W \A\wdtl V=\lambda (\e j\he j) +\jj I \A\jj I-\jj\v (\he j
  \A\jj I).
\ee
\end{Proof}

Similarity transforms as above are used to give bounds on the
maximum eigenvalue. In one case we estimate $\hw$ and $\v$
numerically. For a similarity transform we need to show that
the crucial requirements of $\hw \v=1$ and $\he j\v=1$ are satisfied.
The $j^{\rm th}$ entry in $\v$  can be set to 1, but a numerical
computation yielding $\hw \v=1$ does not necessarily mean $\hw\v=1$.
However, when the numerical computation is done with bounds, one can
show that there are values close to those given so that this sum of
products yields exactly 1.
\begin{Lemma}\label{vwprodlem}
Let
\be
  1-\alpha\leq \sum_0^m v_iw_i\leq 1+\alpha,
\ee
where $0<\alpha<1/2$. With $\lambda:=\alpha/(1-\alpha)$ there exist
$\{w_i^*:\,i=0,\ldots,m\}$ so that
\be
    \sum_0^m v_iw_i^*=1,\quad |w_i-w_i^*|\leq \lambda\, |w_i|,
    \ i=0,\ldots,m.
\ee
\end{Lemma}
\begin{Proof} Let $\sum_{+} v_i w_i$ denote the sum including exactly the terms
with $v_i w_i>0$. Then $\sum_{+} v_i w_i\geq 1-\alpha$.
Define $s_i:=(1-\lambda)w_i$, $t_i:=(1+\lambda)w_i$
for those $i$'s with $v_i w_i>0$; $s_i=t_i=w_i$ otherwise. Now
\be
  \sum_{+} \lambda v_i w_i\geq \alpha\Longrightarrow
  \sum_0^m v_i s_i\leq   \sum_0^m v_i w_i-\alpha\leq 1+\alpha-\alpha=1.  
\ee
Also
\be
  \sum_{+}  v_i t_i\geq \sum_{+} v_i w_i+\alpha\Longrightarrow
    \sum_0^m v_i t_i\geq 1-\alpha+\alpha=1.  
\ee
This implies there exist
$\{w_i^*\colon |s_iw_i|\leq| w_i^*|\leq |t_iw_i|,\,i=0,\ldots, m\}$
so that the sum is unity.
\end{Proof}

%\input car_fp.tex
%&latex
%\iftrue
\iffalse
\nonstopmode
\ifx\pdfoutput\unndefined\else\pdfoutput0\fi
\documentclass[12pt]{article}
\usepackage{amsmath}
\def\dateref{\hfil car-fp.tex 20 February 2016\hfil DRAFT!!}
\textheight=24cm
\textwidth=15.5cm
\def\cD{{\cal D}}\def\cA{{\cal A}}\def\cB{{\cal B}}\def\cS{{\cal S}}
\let\ds\displaystyle
\voffset -2cm
\hoffset -1cm
\newtheorem{Theorem}{Theorem}
\newtheorem{Proposition}{Proposition}
\newtheorem{Lemma}{Lemma}%[section]
\newtheorem{Definition}{Definition}
\newtheorem{Corollary}{Corollary}
\newtheorem{Proof}{Proof}
\def\A{{\bf A}}

\fi
\def\colon{:\,}
\def\frak{\displaystyle\frac}
\def\tA{{\tt A}}\def\tB{{\tt B}}
\def\ta{{\tt a}}\def\tb{{\tt b}}
\def\N{{\bf N}}\def\Z{{\bf Z}}\def\R{{\bf R}}\def\w{{\bf w}}
\def\p{{\bf p}}\def\q{{\bf q}}
\def\e#1{{\bf e}_{#1}}\def\he#1{{\widehat{\bf e}}_{#1}}
\def\tA{{\widetilde {\bf A}}}
\def\wtW{{\widetilde W}}
\let\wdtl\relax
\def\wtV{{\widetilde V}}
\def\crA{{\stackrel{\circ}{\bf A}}}
\def\BB{{\bf B}}
\def\crB{{\stackrel{*}{\bf B}}}
\def\crv{{\stackrel{\circ}{\bf v}}_\lambda}
\def\crw{{\stackrel{\circ}{\bf w}}_\lambda}
\def\crz{{\stackrel{\circ}{\bf z}}_\lambda}
\def\ccz{{\stackrel{\circ}{\bf z}}}
\def\csz{{\stackrel{*}{\bf z}}}
\def\wtV{{\widetilde V}}

\def\v{{\bf v}}
\def\u{{\bf u}}
\def\z{{\bf z}}
\def\jj#1{{#1}^{(j)}}
\def\n{{\bf n}}\def\s{{\bf s}}\def\i{{\bf i}}\def\h{{\bf h}}
\def\ij{{(i,j)}}\def\j{{\bf j}}\def\k{{\bf k}}
\def\hw{{\widehat{\bf w}}}
\def\whi{{\widehat{\bf i}}}
\def\wih{{\widehat{i}}}
\def\whj{{\widehat{\bf j}}}
\def\wjh{{\widehat{j}}}
\let\wh\widehat
\def\hf{{\widehat{f}}}
\def\lrf#1{\|#1\|_r^\infty}
\def\brakt#1#2{\langle\,#1,#2\,\rangle}
\def\D#1#2{D_{#1|#2}}
\def\wM{\widetilde{M}}
\pagestyle{myheadings}\markboth{\dateref}{\dateref}
\parindent 0pt
\def\be{\begin{equation}}\def\ee{\end{equation}}
\def\bea{\begin{eqnarray}}\def\eea{\end{eqnarray}}

\def\prf{\medbreak\noindent{\bf Proof}:\enspace}
\def\qed{\hspace*{\fill}\hbox{\vrule height 7pt \kern-.3pt
    \vbox{\hrule width 7pt \kern6.6pt\hrule width 7pt }\kern-.3pt\vrule height 7pt
    }\par}
\let\vareps\varepsilon
\def\hateps{{\widehat\varepsilon}}
\def\li{\mathop{\rm li}_2\nolimits}

\section{Numerical Approximation}

The behaviour of $f(t)=\ell \log(x)+\sum_1^\infty a_k(t-1)^k$
under the parking and scaling procedure corresponds to the action of
the matrix $\A$ on the coefficients $\{\ell,a_0,a_1,\ldots\}_0$.
We study this behaviour using the largest eigenvalue of $\A$ together
with the right and left eigenvalue pair (\ref{wvdef}). We work with
finite dimensional approximates $\A_m$ of $\A$, and many of the entries of
these are also approximate. Numerical techniques are used to get
estimates of the principal eigenvalue and the corresponding eigenvectors.
We also need error bounds on our approximations to use in connection with
Gershgorin's Theorem to give bounds for the principal eigenvalue of
$\A_m$.

\smallskip 
 
A significant part of the proof depends on numerical calculation
using floating point arithmetic. 
The nonzero real number $x$ has the binary floating point
representation given in terms of a sign bit, an exponent $k$ so that
$2^{k-1}\leq |x|< 2^{k}$ and a $b$-bit integer significand $m$: $2^{-k}|x|$ is
approximated by $m/2^b$, $2^{b-1}\leq m< 2^{b}$. The sign of $x$ is
stored in the sign bit of the floating point representation.
The specific floating point format specifies $b$ and the interval
$[k_{\rm min},k_{\rm max}]$ of allowed values of $k$.
Numbers with absolute value greater than $2^{k_{\rm max}}$ cannot be
represented in this format.
Zero has a special representation; numbers with absolute value less than
$2^{k_{\rm min}}$ are often converted to zero. A nonzero calculation
which yields a result  with absolute value less than $2^{k_{\rm min}}$
is called exponential {\it underflow}; a calculation exceeding the
maximum exponent is called exponential {\it overflow}.
Problems with exponential overflow and underflow do not arise in
the calculations used here, but the limited accuracy of floating point
arithmetic does.

\smallskip

We use $[x]$ to denote the floating point representation of $x$, which
we assume exists in the format specified above. $[|x|]$ is selected so
that $|\,[|x|]-|x|\,|$ is minimal. If $2^{k-1}\leq |x|\leq 2^{k}$,
hence $2^{b-1}\leq 2^{b-k}|x|\leq 2^b$,  there exists an integer $m$,
$2^{b-1}\leq m\leq 2^b$, so that 
\be\label{fpbound}
|2^{b-k}|x|-m|\leq 1/2\Rightarrow|2^{-k}|x|-m2^{-b}|\leq 2^{-b-1}
\Rightarrow \Big|\frac{|x|-[|x|]}{|x|}\Big|\leq 2^{-b},
\ee
since $2^{-k}|x|\geq 1/2$. A special case is when the $m=2^b$, which
will yield the exponent $k+1$ and significand $2^{b-1}$ with the relative error
is close to $2^{-b-1}$. Thus the $b$-bit significand floating point
representation of $x$ has a relative deviation of at most $2^{-b}|x|$
and there may be two distinct floating point representations of $x$
satisfying this relation.

\smallskip

Let $[x],[y]$ be floating point representations and $z:=[x]\bullet [y]$
be the exact result one of the operations $+$, $\times$ or $\div$.
Assume $z$ is defined, nonzero and within the range of the floating
point format.
The value $z$ is not necessarily expressible using $b$-significand bits,
but one can select $[z]$ so that $|z-[z]|/|z|\leq 2^{-b}$. Well
designed floating point routines, e.g. \cite{fpstan}, have this property.
This requires an intermediate result with greater than $b$-bit accuracy.
In the absence of exponent underflow, the result zero arises only if
at least one of the terms in $x\bullet y$ is zero or the operation is
additive with summands equal in magnitude and opposite in sign: if
$x\bullet y$ yields zero, no rounding is needed.

\smallskip

In numerical computations used in the proofs, Intel extended floating
point operations with $b=64$ is used. For basic calculation \textit{round to
nearest} as discussed above is employed. However we also use other rounding
settings to get bounds for the calculation. The setting 
\textit{round up} yields a value greater than or equal to $x\bullet y$;
\textit{round down} yields a value less than or equal to $x\bullet y$.
We use a simplified interval arithmetic with intervals of the form
$x\pm x_r$ and $y\pm y_r$ and compute $z,\,z_r$,
for the operation $x\bullet y$, where $\bullet$ is one of
$+,\times,\div$. We use rounding settings to compute $z_r$ from
$(x,x_r)$  and $(y,y_r)$  so that if
$x-x_r\leq s\leq x+x_r$ and $y-y_r\leq t\leq y+y_r$, then
$z-z_r\leq s\bullet t\leq z+z_r$.  The base value of $x\bullet y$
is computed using round to nearest (in the case of equidistant
nearest values, the even significand is chosen).
We then compute the bound $z_r$ using round up and round down.

The use of interval arithmetic also allows us to bound the maxima and
minima of an arithmetic formula over an interval by dividing the
interval into subintervals and searching, using interval arithmetic
on the subintervals.

\section{Numerical Bounds for the Largest Eigenvalue}

The discussion below uses Gershgorin's Theorem and similarity transforms
involving numerical estimates. It is not sufficient to compute the
transform: we need good bounds on the calculations. Thus the
computations are done with interval arithmetic. One step involves
computing with bounds the transform of the $8\times8$ matrix $\A_7$
which transforms the log term $\ell$ and the polynomial coefficients $a_0,\ldots,a_6$. This is used together
with theoretical bounds involving $\A_m$, $m>7$, to get an interval
containing the maximum eigenvalue.  We obtain bounds for the
right and left eigenvectors for $m>7$. Standard methods are available
for solving linear equations and eigenvalue problems, but usually
these do not provide error estimates. We have combined standard methods
with routines for simplified interval arithmetic to supply the needed
bounds. Our program is written for the Free Pascal
(http://www.freepascal.org/)
compiler for the Intel x86 or x86-64 processor. Our source code is
is available at
\href{http://maths.ucd.ie/renyi}{http://maths.ucd.ie/renyi}.
This  file also contains an estimate of $R_\half$ to
approximately 1000 digits.

The finite dimensional approximation $\A_m$, $m>4$, has a single
eigenvalue $\lambda_m$ of maximum magnitude; other points of the spectrum
have absolute value less than $0.5$. The value of $\lambda_m$
for $m>15$ changes by less than $10^{-10}$.  Computed values of $\lambda_m$
strongly suggest convergence.   Gershgorin's Theorem provides a direct
proof.

\begin{Theorem}
Let $\A=\{a_{ij}\colon i,j=0,\ldots,m\}$ be a matrix with
complex entries. Each eigenvalue of $\A$ is contained in one
of the discs
\be
 \{z\in {\bf C}\colon |z-a_{jj}|\leq\sum_{i\neq j} |a_{ij}|\},
 \ j=0,\ldots,m.
\ee
If the union of the discs about $j_1,\ldots,j_k$ has an empty
intersection with all the remaining discs, this union contains exactly
$k$ eigenvalues, counting multiplicity.
\end{Theorem}
See \cite{Wilk} for details.
\bigskip

In the discussion below, the term \textit{disc} means Gershgorin disc. 
Two matrices related by a similarity transform have the same eigenvalues
and the eigenvectors are related by similarity. We use the above theorem in
connection with some similarity transforms to deduce bounds.
The similarity transforms we use are based on the right and left
eigenvectors associated with the maximum eigenvalue of $\A_m$ for various
$m$. For simplicity of notation the subscript $m$ is often omitted,
when it can be inferred from context.
We write $\v_\lambda$ and $\w_\lambda$ for these eigenvectors. 
We initially require $(\v_\lambda)_0=(\w_\lambda)_0=1$, but for use with
(\ref{wvdef}), $(\v_\lambda)_0=1$, but $\w_\lambda$ is scaled
so that $\widehat\w_\lambda \v_\lambda=1$. Define $\crA$ with
\be\label{crAdef}
\crA_{ij}= \A_{ij},\; i,j=1,2,3,\ldots\,,
\ee
$\crA$ is the matrix $\A$ with row and column $0$ omitted. We denote
the associated truncated eigenvectors with the 0 entries removed
by $\crv$ and $\crw$. The eigenvalue equations
$\A\v_\lambda=\lambda\v_\lambda$, $\widehat\w_\lambda \A=\lambda\widehat\w_\lambda$
yield the linear equations
\be\label{wveqns}
(\lambda-\crA)\crv=\crA_{\bullet0};\quad
\hcrw(\lambda-\crA)=\crA_{0\bullet},
\ee
where $\crA_{\bullet0}$ and $\crA_{0\bullet}$ denote the column and row of
$\A$ indexed by 0 with the entries indexed by zero removed.
We use the $\ell^1$ norm on column vectors and
the $\ell^\infty$ norm on row vectors: 
\be
\|\v\|=\sum_i|\v_i|\ ;\qquad\|\w\|_\infty=\max_i\{|\w_i|\};
\ee
\be
  \|\v_\lambda\|=\|\crv\|+1;\qquad\|\w_\lambda\|_\infty=\max\{1,\|\crw\|_\infty\}.
\ee
The associated matrix norm, $\|\crA\|$, is the maximum $\ell^1$ norm of
the columns of $\crA$. Note that $\crA$ is upper triangular.
The absolute sum of column $k+1$  is bounded by
\be\label{colsumbnd}
  \frac1k+2^{-k-1}(1+\frac1k)\sum_{j=0}^k\binom{k}{j}\leq 1
\ee
for $k\geq 3$. Checking the columns for $k<3$
shows that the supremum of the absolute column sums of $\crA$
occurs at the first column: $\|\crA\|=(\A)_{1,1}=0.5+\log 2$.
For $\lambda>(\A)_{1,1}$, which is the case for the maximum eigenvalue,
we have
\be\label{pinverse}
 ( \lambda-\crA)^{-1}=
  \sum_0^\infty\lambda^{-k-1} (\crA)^k ,\quad
  \|(\lambda-\crA)^{-1}\|\leq \frac 1{\lambda-\A_{1,1}}.
\ee
In the case of $\A_m$, $\crA$ is indexed
by $1\ldots m\times 1\ldots m$ and $\lambda=\lambda_m$, the maximum
eigenvalue of $\A_m$.
Because $\lambda-\crA$ is non-singular, $\A\v=\lambda \v$ with
$\v_0=0$ implies $\v=0$.  Thus, if $\A\v_\lambda=\lambda \v_\lambda$
with $\v_\lambda\ne0$, we can scale $\v_\lambda$ to make
$(\v_\lambda)_0=1$. 
The inverse of $(\lambda-\crA)$ can be computed by the power series above
or by using the upper triangular property of $\crA$. One can use back
substitution to solve for $\crv$. The upper triangular forward
substitution for $\crw$ depends on $\lambda$ but the evaluation of
$(\crw)_k$ depends only on $\lambda$ and $\A_{i,j}\colon i,j\leq k$.

\begin{Lemma}  Let $\sigma_i=i$ for $i>1$, $\sigma_1=0$,
$\sigma_0=1$. Then $(-1)^{1+\sigma_i}(\crv)_i>0$ for each integer
$i\geq 0$.
\end{Lemma}

\begin{Proof} Define the matrix $\BB$ by
\be
   \BB_{ij}:=(-1)^{\sigma_i-\sigma_j}\A_{ij}.
\ee
Then all the $\BB$ entries $\BB_{ij}$ for $i\neq 1$ except
$\BB_{2,2}$ are nonnegative.
Define $\crB$ as $\crA$ in (\ref{crAdef}) except that row and column 1
are removed. Then all entries of $\crB$ are nonnegative except the
term coming from $\BB_{2,2}=-1/4$.
Let $\BB\z=\lambda\z$ with
$\z_1=1$. Let $\csz$ be $\z$ with the 1 entry removed.
$\crB_{\bullet1}$ has a single nonzero entry $\crB_{0,1}=1$.
As in (\ref{colsumbnd}) one can show that $\|\crB\|<1$, then
\be
   \csz=\sum_0^\infty(\lambda+1)^{-k-1}(\crB+I)^k (\crB_{\bullet1}),
\ee
so the entries of $\csz$ are strictly positive including $\z_1=1$.
We  have $\v_0=1$ by definition so
$\v_j=(-1)^{\sigma_j-\sigma_0}z_j/z_0$
with $\sigma_0=1$ from
the similarity transform relating $\A$ and $\BB$, which uses the
self-inverse diagonal matrix with diagonal entries $(-1)^{\sigma_j}$.
\end{Proof}

\begin{Lemma} For $m\geq 7$, if $\lambda_m$ satisfies 
$1.232\leq\lambda_m\leq 1.234 $, then for
$\A_m \v_{\lambda_m}=\lambda_m\v_{\lambda_m}$ with $\v_0=1$
\be\label{sumandfirst}
  2.03<\sum_{i=0}^m|\v_i|<2.08,\quad  -0.76<\v_1<-0.70.
\ee
\end{Lemma}

\begin{Proof}
We consider $\A_m$ for some $m\geq 7$.
Let $\z$ be a row vector with $|(\z)_i|\leq 1$. Then
\be
  \z\v=(\z)_0(\v)_0+\ccz\crv
  =(\z)_0(\v)_0+\Big(\ccz(\lambda-\crA)^{-1}\Big)\crA_{\bullet0}.
\ee
With $(\crA_{\bullet0}^{(7)})_{i}=\A_{i}$ for $1\leq j\leq 7$ and
$(\crA_{\bullet0})_{i}=0$ for $i>7$, we have
\be
\left|\Big(\ccz(\lambda-\crA)^{-1}\Big)\crA_{\bullet0}- 
\Big(\ccz(\lambda-\crA)^{-1}\Big)\crA_{\bullet0}^{(7)}\right|
\leq  \| (\lambda-\crA)^{-1} \|
|\crA_{\bullet0}-\crA_{\bullet0}^{(7)}|
\ee
First we select $\z$ with $\z_0=1$, $\z_1=-1$ and for $i>1$ $(\z)_i=(-1)^{i+1}$ so
that $|(\z)_i|=1$ and $(\ccz)_i(\crv)_i> 0$
because the signs of $(\crv)_i$ alternate.
We compute using interval arithmetic that
2.032 and 2.072 are bounds for the minimum and maximum of
$\Big(\ccz(\lambda-\crA)^{-1}\Big)\crA_{\bullet0}^{(7)}$ over
$1.232\leq\lambda\leq 1.234$. Next from (\ref{col0bnd})
\be\label{corr}
  | \|\v\|-1-\Big(\ccz (\lambda-\crA)^{-1}\Big)\,\crA_{\bullet0}^{(7)}|
  \leq\frac{\| \crA_{\bullet0}-\crA_{\bullet0}^{(7)}\|}
  {1.232-0.5-\log(2)}<\frac {2.25\times 3^{-7}}{49\times0.03885}
  <0.00055,
\ee
from which follow the first pair of inequalities in (\ref{sumandfirst}).
For the second pair of inequalities we set  $\z_i=0$ for
$i \neq 1 $, $\z_1=1$. Then as above
\be\label{corrr}
  | \v_1 -\Big(\ccz (\lambda-\crA)^{-1}\Big)\,\crA_{\bullet0}^{(7)}|
  <0.00055.
\ee
We compute using interval arithmetic the values -0.711 and -0.751
as  upper and lower bounds for
$\Big(\ccz (\lambda-\crA)^{-1}\Big)\,\crA_{\bullet0}^{(7)}$ 
over $1.232\leq\lambda\leq 1.234$.  Following this, 
(\ref{corrr}) shows that the second pair of inequalities obtain.
\end{Proof}

\begin{Lemma} For $m\geq 7$, if $\lambda_m$ satisfies 
$1.232\leq\lambda_m\leq 1.234 $, then for
$\A_m \v_{\lambda_m}=\lambda_m\v_{\lambda_m}$ with $\v_0=1$ and
$\w_{\lambda_m}\A_m =\lambda_m\w_{\lambda_m}$ with $\w_0=1$
\be\label{prodnrmbnd}
  19.2<\sum_{i=0}^m \w_i \v_i<21.4,
  \quad 1.14\leq \frac{\max_{0\leq i\leq m}\{|\w_i|\}}
  {\sum_{i=0}^m \w_i \v_i}\leq 1.35.
\ee
The matrices $W,\,V$ defined by (\ref{wvdef}) satisfy
\be\label{wvnorms}
  \|W\|\leq 1+\|\w\|_\infty\leq 2.35,\quad
  \|V\|\leq 1+\|\v\|\|\w\|_\infty\leq 3.81.
\ee
\end{Lemma}
\begin{Proof}
We have
\be
 \sum_{i=0}^m \w_i \v_i=1\times 1+
 \crw (\lambda-\crA)^{-1}\crA_{\bullet0}.
\ee
\be
  |\sum_{i=0}^m \w_i \v_i-1-\crw (\lambda-\crA)^{-1}\crA_{\bullet0}^{(7)}|
  \leq 1*\|(\lambda-\crA)^{-1}\|^2
  \,\|\crA_{\bullet0}-\crA_{\bullet0}^{(7)}\|,
\ee
since in the $\ell^\infty$ norm $\|\crA_{0\bullet}\|=1$ and
  $\crw=\crA_{0\bullet} (\lambda-\crA)^{-1}$.
We use interval arithmetic to find 
\be
19.3<\sum_{i=0}^7 \w_i \v_i<21.3.
\ee
Then $2.25\times 3^{-7}/(49(1.232-0.5-\log(2))^2)<0.014$
gives a bound for the terms greater than 7. 
For the bound on $\w$, note that with $(\w_\lambda)_0=1$ forward
substitution in (\ref{wveqns}) yields 
$(\w_\lambda)_1=-1/(\lambda-A_{1,1})$ and (\ref{pinverse}) then shows
that $\|\w_\lambda\|=|(\w_\lambda)_1|$. 
With 19.2 and 21.4 lower and upper bounds for $\sum\w_i\v_i$
and $1/(1.234-0.5-\log2)$, $1/(1.232-0.5-\log2)$
lower and upper bounds for  
$\max_{0\leq i\leq m}\{|\w_i|\}$, which occurs at $i=1$,
we have the last pair of inequalities of (\ref{prodnrmbnd}).
Now $W=\e j\hw-\jj\v\he j+\jj I$ with $j=0$. The norm of column
$0$ of $W$  is $\|\v\|-1+1/\sum \w_i\v_i< 2.08$. The column with
greatest absolute sum is column 1, so we get the scaled $|\w_1|$ plus 1.
Then \eqref{prodnrmbnd} implies $\|W\|\leq 2.35$. We have
$V=\v(\he j-\jj\hw )+\jj I$ again with $j=0$. The maximum column
absolute sum occurs at $i=1$. It is bounded by $\|\v\|\,\|\w\|+1$.
\end{Proof}

\begin{Proposition}\label{lamest} For $m\geq 7$, $\A_m$ has a single real eigenvalue in the
interval $1.232891\pm0.0007$. All remaining eigenvalues have
absolute values less than 1.002.
\end{Proposition}

\begin{Proof}
First we consider $\A_7$. We compute
approximate right eigenvector $\v$ and
left eigenvector $\hw$ for $\lambda\approx 1.232891$ with $\w_0=1$
and $\hw\,\v\approx 1$.
First we compute with interval arithmetic that
$|\sum \w_i\v_i-1|<7\times 2^{-64}$, so we can apply Lemma \ref{vwprodlem} to
conclude that there exists $\w^*$ so that $\sum \w_i \v_i^*=1$ and
$|\w_i-\w_i^*|\leq \vareps|\w_i|$ with $\vareps=8\times 2^{-64}$. We define
$W$ and $V$ by (\ref{wvdef}). Using
interval arithmetic we deduce a disc about 1.232891 with
radius less than $53\times 2^{-64}$. The absolute sum of each of the
remaining columns is bounded by $1+46\times 2^{-64}$,
which implies that if the complex number $z$ is in any  discs except the one about
1.232891, then $|z|\leq 1.001$. Thus there is a single eigenvalue with
magnitude greater that 1.001. This must be real; otherwise we would
have a conjugate pair of eigenvalues of magnitude greater than 1.001.

Next we consider the case of $\A_m$, $m>7$. Extend $W$ and $V$ to $\A_m$
as follows using $h=7$:
\be\label{wtWVdef}
  \wtW_{ij}=
  \begin{cases}W_{ij}& \text{if $0\leq i,j\leq h$;}\\
                  0& \text{if $ i,j> h,\,i\neq j$};\\
		  8& \text{if $i=j> h$;}
	     \end{cases}
\quad	  
   \wtV_{ij}=\begin{cases}V_{ij}& \text{if $0\leq i,j\leq h$;}\cr
		                0& \text{if $ i\neq j$.}\cr
		               \frac18& \text{if $ i=j> h$}.\cr
	      \end{cases}
\ee
We have the similarity transform $\wtW \A\wtV$. 
We need bounds for the discs about
$(\wtW \A\wtV)_{ii}$ for $i>7$ and bounds for the contributions of
the terms $(\wtW \A\wtV)_{ij}$ to the discs about 
$(\wtW \A\wtV)_{ii}$ for $i\leq 7$. For $j>7$
\be\label{88125}
  \sum_{i=0}^m |(\wtW \A\wtV)_{ij}|\leq
  \sum_{i=0}^6 \frac{\|W\|}{8} \,|\A_{ij}|
  +\sum_{i=7}^j |\A_{ij}|
    <\frac{2.35}{8} (0.5+\frac26)+0.5+\frac16<1.
\ee
For $j\leq 7\leq i$, $ \A_{ij}\neq 0$ only if $j=0$. Thus for $j\leq 7$, 
\be\label{00064}
\sum_{i=7}^m  |(\wtW \A\wtV)_{ij}|\leq
8\|V\|\sum_{i=7}^\infty \A_{i0}<8*3.81*2.25*3^{-7}/7^2<0.00064,
\ee
using (\ref{col0bnd}) and (\ref{wvnorms}).
For $\wtW \A\wtV$ the  disc about  1.23289 has radius
less than $0.00065$; any $z$ belonging to the remaining discs has
$|z|<1.001+0.00065$.
\end{Proof}

\begin{Proposition} \label{limpro}
Let $(\lambda_n,\v_n,\w_n)$ denote the maximum eigenvalue and
corresponding right and left eigenvectors of $\A_n$, where the
$(\v_n)_1=(\w_n)_1=1$.  Then the sequence $\{\lambda_n\}$ converges,
while $\{\v_n\}$ converges in $\ell^1_\rho$, $1\leq \rho\leq 3$ and 
$\{\w_n\}$ converges in $\ell^\infty$.
\end{Proposition}
\begin{Proof}
For $\A_m$, $m>7$ and $n>m$, define $W_m$ and $V_m$ using  $\zeta\w$ for $\w$ in
(\ref{wvdef}) with $\zeta=(\w\v)^{-1}$. By (\ref{pinverse}) and the
upper triangular form of $\crA$ with $\|\crA_{\bullet0}\|_\infty=1$,
\be
\|\w\|_\infty\leq 1/( \lambda_m-\A_{1,1}),
  \w_1=1/( \lambda_m-\A_{1,1}),
  \|\zeta\w\|_\infty= ( \lambda_m-\A_{1,1})^{-1}(\w\v)^{-1}.
\ee
From (\ref{sumandfirst}) have $\v_0=1$, $-0.76<\v_1<-0.70$, $\sum_i |\v_i|<2.08$,
so $|\v_i|\leq 1$ for $i=0,\ldots,m$.
Then for $n>m$ define $\wtV$ and $\wtW$ using (\ref{wtWVdef})
with $h=m$. First we consider  $(\wtW \A_n\wtV)_{ij}$ for
$0\leq i,j\leq m$. These entries are given by (\ref{eq:47}). For
$0\leq i,j\leq m$, $i\times j=0$ the values are zero except the (0,0) entry
is $\lambda_m$. For $1\leq i\leq m$, $1\leq j\leq m$,
\be
(\wtW \A_n\wtV)_{ij} =(\A_n)_{ij}-\v_i(\A_{0,j}).
\ee
Note $(\A_{0,j})=0$ for $j>1$, so the entries equal $(\A_n)_{ij}$
when $1\leq i\leq m$ and $2\leq j\leq m$. The argument for (\ref{colsumbnd})
and estimates of the form (\ref{00064}) show that
\be
2\leq j\Rightarrow\sum_{i=0}^n(\wtW \A_n\wtV)_{ij}< 1.002.
\ee
Only the column indexed by $j=1$ is different. 
\be
  (\wtW \A_n\wtV)_{0,1}=0,\ (\wtW \A_n\wtV)_{1,1}=0.5+\log 2+\v_1,
  (\wtW \A_n\wtV)_{i,1}=\v_i,\,1<i\leq m.
\ee
Since $\sum_0^m\|\v_i\|<2.08$ and $-0.76<\v_1<-0.70$,
\be
  \sum_{i=0}^m|(\wtW \A_n\wtV)_{i1}|=\sum_0^m\|\v_i\|-1+0.5\log2
  -2\v_1<0.874.
\ee
The argument of (\ref{00064}) implies
$\sum_{i=0}^n|(\wtW \A_n\wtV)_{i1}|<0.875$.
For $j\leq m$, 
\be
\sum_{i=m+1}^n |(\wtW \A_n\wtV)_{ij}|
\leq 8\|V\|\frac{2.25\times 3^{-m}}{m^2}.
\ee
The above is a bound for the radius of the  disc about
$\lambda_m$. For $n>m$, $\lambda_n$ is inside this disc. Thus shows that
for $n,k>m$ using (\ref{wvnorms}),
\be
  |\lambda_n-\lambda_k|\leq 16 \|V\|\frac{2.25\times 3^{-m}}{m^2},
\ee
hence $\{\lambda_n\}$ is a Cauchy sequence. Then (\ref{wveqns})
and (\ref{pinverse}) show that $\crv$ and $\crw$ converge as $\lambda_n$
converges, hence $\{\v_{\lambda_n}\}$ converges in $\ell^1$
and $\{\w_{\lambda_n}\}$ converges in $\ell^\infty$. To see that
$\{\v_{\lambda_n}\}$ converges in $\ell^1_\rho$, note that the
transformed matrix entries $(\A)_{jk}\rho^{j-k}$ in $\crA$ are not greater
than those in $\crA$ for $\rho\geq 1$ because $\crA$ is upper triangular.
In (\ref{wveqns}) it follows from (\ref{col0bnd}) that
the transform of $\crA_{\bullet0}$ is in
$\ell_\rho^1$ for $1\leq \rho\leq 3$, so $\crv$ converges in $\ell_\rho^1$.

\end{Proof}

\begin{Remark}
 Proposition~\ref{lamest}  implies the limit $\lambda=\lim \lambda_n$
 lies in the interval $1.232891\pm 0.0007$.
\end{Remark}

\iffalse
\begin{Theorem} Let $\lambda$, $\w_\lambda$ and $\v_\lambda$ be the limits
as above with $\w\v=\zeta$.
Then if $\u\in\ell^1$ satisfies $\w\u\neq 0$, then
\be
\lim_{n\to\infty} \lambda^{-n} \A^n\u=(\w\u)\v/\zeta.
\ee
Also if $\w\u\neq 0$,
\be\label{limrat}
\lim_{n\to\infty} \frac{ \A^n\u}{\| \A^n\u\|}
=\frac{(\w\u)\v}{\|(\w\u) \v\|}.
\ee
\end{Theorem}
\begin{Proof}
Use (\ref{wvdef})  to define $W$ and $V$ on
$\ell^1$. Define $\z$ and $B=(B_{ij})$ by
\be
\z=W\u-\frac1\zeta(\w\u)\he0,\quad B=
\begin{cases}
     (W \A V)_{ij} & \text{if $ij\neq 0$};\cr
      0 & \text{if $i=0$ or $j=0$}.\cr
\end{cases}					       
\ee
Then
\be
   (W \A V)^n\u=\lambda^n \frac1\zeta(\w\u)\he0 +\lambda^n B^n\z.
\ee
The estimates in Proposition \ref{limpro} show that the absolute sum of
each column of $B$ is not greater than 1, i.e., $\|B\|\leq 1$. Then
\be
 \lambda^{-n}(W \A V)^n \z=
 \lambda^{-n}(W \A V)^n W\u- \frac1\zeta(\w\u)\he0
  =\lambda^{-n}B^n\z.
\ee
Now $\A^n \u=V(W \A V)^n W\u$ and $\|B\|\leq 1$ and $V\he0=\v$ so
\be
  \|\lambda^{-n} A^n\u-\frac1\zeta (\w\u)\v\|\leq \lambda^{-n}\|V\|\,\|z\|
  \leq  \lambda^{-n}\|V\|\,(\|W\|+\frac{\|w\|}{\zeta})\|u\|.
\ee
For (\ref{limrat})
$\|\lim \lambda^{-n} \A^n\u\|=\lim \|\lambda^{-n} \A^n\u\|$; one
takes the limit of the ratios.
\end{Proof}
\fi

\begin{Theorem} Let $\lambda$, $\w_\lambda$ and $\v_\lambda$ be the limits
as above with $\w\v=\zeta$.
Then if $\u\in\ell^1$ satisfies $\w\u\neq 0$, then
\be
\lim_{n\to\infty} \lambda^{-n} \A^n\u=(\w\u)\v/\zeta.
\ee
Also if $\w\u\neq 0$,
\be\label{limrat}
\lim_{n\to\infty} \frac{ \A^n\u}{\| \A^n\u\|}
=\frac{(\w\u)\v}{\|(\w\u) \v\|}.
\ee
\end{Theorem}
\begin{Proof}
Use (\ref{wvdef})  to define $W$ and $V$ on
$\ell^1$, where one uses $\w/\zeta$ in the expressions for $V$ and $W$
to correspond to the $\w$ used here. 
Write $\u=\v(\w\u)/{\zeta}+\u-\v(\w\u)/\zeta$ and define
$\z$ and $B=(B_{ij})$ by 
\be
\z=W\u-\frac1\zeta(\w\u)\e0,\quad B=
\begin{cases}
     (W \A V)_{ij} & \text{if $ij\neq 0$};\cr
      0 & \text{if $i=0$ or $j=0$}.\cr
\end{cases}                                            
\ee
Then $W\u=(\w\u/\zeta)\e0 +\z$ with $\he0 \z=0$, so
\be
   (W \A V)^n W\u=\lambda^n \frac1\zeta(\w\u)\e0 + B^n\z.
\ee
The estimates in Proposition \ref{limpro} show that the absolute sum of
each column of $B$ is not greater than 1, i.e., $\|B\|\leq 1$. Then
\be
\lambda^{-n}\A^n \u=
V\lambda^{-n}(W \A V)^n W\u= (\w\u/\zeta)V\e0+V\,\lambda^{-n}B^n\z.
\ee
Now $\|B\|\leq 1$ and $V\e0=\v$ so
\be
  \|\lambda^{-n} \A^n\u-\frac1\zeta (\w\u)\v\|\leq
  \lambda^{-n}\|V\|\,\| B\|^n \|\z\|
  \leq  \lambda^{-n}\|V\|\,(\|W\|+\frac{\|\w\|}{\zeta})\|\u\|.
\ee
For (\ref{limrat})
$\|\lim \lambda^{-n} \A^n\u\|=\lim \|\lambda^{-n} \A^n\u\|$; one
takes the limit of the ratios.
\end{Proof}

In this section we expressed the transformation of
$\ell \log(t)+\sum_0^\infty a_k (t-1)^k$ as the matrix operation $\A$.
Now we consider the reverse process. Let
$(\ell,a_0,a_1,\ldots)\in\ell^1$ and define
\be
  g(t)=\ell^*\log(t)+\sum_{k=0}^\infty a_k^*(t-1)^k.
\ee
We write $\A g$ to denote the corresponding function from the
coefficients $\A(\ell,a_0,a_1,\ldots)$. The definition of
$\A$ implies
\be
  g\geq 0\Rightarrow \A g \geq 0,\quad
  \int_0^2 \A g\,dt=\int_0^1 g\,dt=2\int_1^2  g\,dt.
\ee
In particular, if $\A g=\lambda g$ and $\int_0^2 g\,dt\neq 0$,
\be
  \lambda\int_0^2 g\,dt=\int_0^2 g\,dt+\int_1^2  g\,dt
  \Rightarrow\lambda=1+\int_1^2  g\,dt\Big/\int_0^2  g\,dt.
\ee

Here we show that the appropriately scaled function corresponding to the
eigenvector $\v$ satisfies the hypotheses of Proposition
\ref{prop:props}.

\begin{Corollary}\label{fstarcor}
There exists a probability density $f^*$ on $(0,2]$ of the form
\be\label{fstardef}
   f^*(t)=\ell^*\log(t)+\sum_{k=0}^\infty a_k^*(t-1)^k
\ee
so that $f^*$ satisfies (\ref{eq:int_ss}) with $1+C=\lambda$. The coefficients
$(\ell^*,a_0^*,a_1^*,\ldots)$ are in $\ell^1_\rho$ for $1\leq \rho\leq 3$.
Starting from any probability density of the form
\be
   f_0(t)=\ell\log(t)+\sum_{k=0}^\infty a_k(t-1)^k
\ee
with $f_s(t)$ given by (\ref{eq:int_ev}) for each $0<y<2$ we have
\be
  \lim_{s\to\infty}\sup_{y\leq t\leq2} |f_s(t)-f^*(t)|=0.
\ee
\end{Corollary}
\begin{Proof}
  First let $f_0=0.5$ on $[0,2]$ corresponding to  $\u$ with $\ell=0$,
$a_0=0.5$ and $a_k=0$ otherwise. We know that $\A^n f_0 \geq 0$ and
$\lim_n \lambda^{-n} \A^n\u=(\w\u)\v/\zeta$.
Define $g(t)$ by
\be
  g(t)=-\v_0\log(t)-\sum_{k=0}^\infty \v_{k+1}(t-1)^k.
\ee
Then $\A g=\lambda g$. In the representation as functions
\be
\lim_{n\to\infty} \lambda^{-n} \A^n f=-\frac{\w_1}{2\zeta}g,
\ee
because convergence of the coefficients implies uniform convergence on
compact intervals of $(0,2]$. From (\ref{prodnrmbnd}) and $\w_1=-\|\w\|$
we have $-\w_1=1/(\zeta(\lambda-0.5-\log(2)))\geq 1.14$,
so the limit of $\lambda^{-n}\A^n f$ is $c\,g(t)$
with $c\geq 0.57$, hence $g(t)\geq 0$.
The $\lim_s f_s$ is the multiple of $g$ with
$\int_0^2 g\,dt=1$:   $f^*(t)=g(t)/\int_0^2 g\,dt$.
Now for a general initial probability density of the given form,
$\A^k f_0$ will be strictly positive for some $k>0$
$f_s=\widehat f_s+\vareps$. Then
$\lambda^{-k} A^k (\widehat f_s+\vareps)\geq \lambda^{-k} A^k \vareps$. Then the
$\lim_n\lambda^{-n}\A^n f_0$ cannot be the zero function.
By taking ratios we find that $\lim_{s\to\infty}f_s=f^*$.
\end{Proof}

\def\hT{{\widehat T}}
\def\tT{{\widetilde T}}
\def\cM{{\cal M}}\def\cMp{{\cal M}^{+}}\def\cMo{{\cal M}_1^{+}}

\begin{Remark}
  By virtue of \eqref{decay}, we now have an estimate of the rate of
  decay of remaining space in the interval with error estimate,
  $R_\half=0.616445\pm0.0035$.  
% Source code used in determining the
%   above estimates
% %, including a calculation of $R\half$ to 1000 digits
% %  using the GNU multiple precision arithmetic library \cite{gmp}
%   is available at
%   \href{http://maths.ucd.ie/renyi}{http://maths.ucd.ie/renyi}.
%   This  file also contains an estimate of $R_\half$ to
%   approximately 1000 digits.

\end{Remark}

\section{The Scaled Rényi Parking Transform}\label{sec:gen}

We have shown in Corollary~\ref{fstarcor} that a probability density of the form \eqref{eq:fform}
with coefficient vector $\{a_k\} \in \ell^1$ converges under
the iteration to a density $f^*$.    Here we consider more general
 distributions on $[0,2]$.   
% Note that on the closed
% interval $[0,2]$, the point mass at $0$ is a fixed point of the
% iteration, so does not converge to $f^*dx$.  
In this section, we express the iteration on non-negative
measures.  This iteration  has $\delta_0$ as a fixed
point,  but,  in the case of probability
measures on $(0,2]$, is equivalent to that for
cumulative distribution functions given in \eqref{eq:distev}.  In
particular, we show that the point measure $\delta_x$ for $x\in (0,2]$ converges to
$f^* dx$.

\begin{Definition}
$C[0,2]$ denotes the continuous functions on $[0,2]$; $\cM$ denotes
the bounded Borel measures on $[0,2]$;
$\cMp$ denotes the non-negative  measures in $\cM$;
$\cMo$ denotes the measures $\mu\in\cMp$ for which $\mu[0,2]=1$.
The integral of $g\in C[0,2]$ with respect to $\mu\in\cMp$ is denoted by
\be
  \brakt g \mu\equiv \int_0^2 g(x)\mu(dx)\equiv \int g \,d\mu.
\ee
For $\mu,\nu\in\cM$ we say that $\mu\leq\nu$ if $\nu-\mu\in\cMp$.
\end{Definition}
The value of $\brakt g \mu$ for all $g\in C[0,2]$ determines $\mu\in\cMp$ uniquely.
For bounded Borel functions $f,g:[0,2]\to\R$, we have the measures
$f(x)\,dx\leq g(x)\,dx$ if and only if $f\leq g$ a.e.

We define the operator $U\colon\cMp\to\cMp$ by
\be
  \brakt g {U\mu}=\int_0^1 g(2x)\mu(dx). %%%\brakt{g(2x)}\mu.
\ee
Thus $U\mu$ takes the restriction of $\mu$ to $[0,1]$ and scales it to
an element of $\cMp$ using the operator $x\mapsto 2x$.

\bigskip

For each $x\in (1,2]$ we define the measure $\nu_x\in\cMp$ by
\be\label{nuxdef}
  \brakt g {\nu_x}=\frac 1{x-1}\int_0^2 g(y)  I_{[0,2(x-1)]}(y)\,dy
  = \frac 1{x-1}\int_0^{2(x-1)}g(y)\,dy,
\ee
so $\nu_x$ is Lebesgue measure on $[0,2(x-1)]$ scaled by $1/(x-1)$. Note
$\nu_x[0,2]=2$. Here we use the indicator function
\be
  I_{[0,2(x-1)]}(y)=I_{y\leq 2(x-1)}=I_{x\geq 1+y/2}=
  \begin{cases} 1 & \mbox{if } 0\leq y\leq 2(x-1)\leq 2;\cr
                0 &\mbox{otherwise}.\cr
  \end{cases}		
\ee
Given $\mu\in\cMp$, $g\in C[0,2]$ we define $V\mu$ by
\be
  \brakt g {V\mu}:=\int_{1^+}^2 \brakt g{\nu_x}\mu(dx)
  =\int_{1^+}^2 \Big(\int_0^{2(x-1)} \frac{g(y)}{x-1}\,dy\Big)\,\mu(dx).
\ee
Note that the integral is on $(1,2]$; the value $\mu\{1\}$ is not
included. The treatment of $\mu\{1\}$ is significant. The measure concentrated
at $1$, $\delta_1$, has $V \delta_1=\delta_2$. If we included the point
$1$ in the definition of $V$, then the image should be $2\delta_0$.
We define
\be
T=U+V
\ee
as the basic unscaled  Rényi transformation. The following is elementary.
\begin{Proposition}
  If $\mu\leq\nu$ then $U\mu\leq U\nu$, $V\mu\leq V\nu$ and $T\mu\leq T\nu$.
\end{Proposition}

The point 1 is a discontinuity in $T$: for $x\leq 1$, $T\delta_x$
is $\delta_{2x}$; for $x>1$, $T\delta_x=\nu_x$. % (\ref{nuxdef}).
Note that $\delta_0$ is a fixed point of $T$.

\bigskip

We define the renormalised transform
for $\mu\in \cMp$, $\mu\neq 0$,
\be\label{areascaled}
  \hT\mu= \frac1{(T\mu)[0,2]}T\mu,
\ee
so that $\hT\mu\in \cMo$. The transformations, $T,U,V$ are linear. The
transformation $\hT$, called the \textit{area scaled}  Rényi
transformation, is nonlinear, but $\hT$ maps $\cMo$ into itself. Note that
$\hT^n \mu =\hT(T^{n-1}\mu)$, so one can work with $T$ and normalise
in the last step. For the  probability density $f^*$ from
Corollary~\ref{fstarcor}, 
with $C^*=\int_1^2 f^*(x)dx$, we define
\be\label{spectrumscaled}
  \tT\mu= \frac1{1+C^*}T\mu.
\ee
This just a different normalization of $T$ which leaves the measure
$f^*dx$ fixed. 
Notice $\hT f^*\,dx=f^*dx$ as well. 
Also $\mu\leq\nu\Rightarrow \tT\mu\leq \tT\nu$, but this need not hold
for $\hT$.
\bigskip

Now we consider the case when $\mu\in\cMp$ has  a density
$\mu(dx)=f(x)dx$. For $g\in C[0,2]$
\be
  \brakt{g}{U \mu}=\int_0^1 g(2y) f(y)\, dy=
  \int_0^2 g(x) f(x/2)\, dx/2,
\ee
so $U$ acts on $f$  by $f\mapsto \frac12 f(x/2)$ for $x\in[0,2]$.
Next consider $V$.
\be
 \brakt{g}{V \mu}=
\int_{1^+}^2 \Big(\int_0^{2(x-1)} \frac{g(y)}{x-1}\,dy\Big)f(x)\,dx=
\int_{1^+}^2\int_0^2 I_{y\leq 2(x-1)} \frac{g(y)}{x-1} f(x)\, dy\, dx.  
\ee
Changing the order of integration we get
\be
 \brakt{g}{V \mu}=
\int_0^2\int_{1^+}^2 I_{x\geq 1+y/2} \frac{g(y)}{x-1} f(x)\, dx\, dy=
\int_0^2 g(y)\int_{1+y/2}^2 \frac {f(x)}{x-1}\,dx\,dy\,. 
\ee
Thus the action of $V$ on the density $f$ is given by
\be
  f\mapsto \int_{1+x/2}^2 \frac {f(y)}{y-1}\,dx
\ee
Compare the above with (\ref{eq:int_ev}).
The following shows that for any $\mu\in\cMp$ with $\mu(0,2]>0$,
$T^n\mu-\vareps\,dx\in\cMp$ for some $n\in\N$ and $\vareps>0$.
\begin{Proposition}\label{prop:mp}
Let $\mu\in\cMp$ satisfy $\mu(0,2]>0$. Then there exists $n$,
$\vareps>0$ and $\nu\in\cMp$ so that $T^n\mu= \nu +\vareps dx$ and
$\alpha>0$ so that  $T^{n+1}\mu\geq \alpha \,f^*dx$.
\end{Proposition}
\begin{Proof} Since $\mu(0,2]>0$, we can iterate $T$ as necessary to get
$(T^j\mu)(1,2]>0$.
Given $T^j\mu(1,2]>0$, there exists $k$ so that
$(T^j\mu)[1+2^{-k},2]>0$. $V T^j \mu$ integrated over
$[1+2^{-k},2]$ is an integral of $\nu_x$ over $[1+2^{-k},2]$.
The minimal density of $\nu_x$ for $1+2^{-k}\leq x\leq 2$ occurs at $1$
and equals 1. Thus 
$T^{j+1} \mu \geq (T^j\mu)[1+2^{-k},2]I_{[0,2^{-k+1}]}dx$.
We then have $UT^{j+k} \mu\geq \vareps I_{[0,2]}dx$
with $\vareps=(T^j\mu)[1+2^{-k},2]$. Note
\be
  V T^{j+k}\mu\geq \int_{1+x/2}^2\frac
  {\vareps\,dy}{y-1}dx=\vareps\log\frac2x,
\ee
and $f^*$ is a linear combination of $\log x $ and a bounded function,
so there exists $\alpha>0$ so that
\be
T^{j+k+1} \mu \geq \vareps(1+\log\frac2x)dx\geq \alpha f^*dx.
\ee
\end{Proof}

\begin{Theorem} Let $f(x)=D\log(x)+g(x)$ be a probability density on
$(0,2]$ where $D$ is a constant and $g(x)$ is continuous on
$[0,2]$. Then $\lim_{n\to\infty}\hT^n f\,dx=f^*\,dx$.
\end{Theorem}
\begin{Proof} The continuity of $f$ on $(0,2]$ implies via
\eqref{eq:int_ev} the continuity of $\hT^n f$.
Since $g$ is continuous, for $\vareps>0$ there exists a polynomial
$q_\vareps(x)$ so that
\be
  \sup_{0<x\leq 2} |f(x)-D\log(x)-q_\vareps(x)|\leq\vareps
\ee
because the polynomials are dense in the set of continuous functions
on $[0,2]$. We have
\be
\tT^n (D\log(x)+q_\vareps(x)-\vareps)\leq\tT^nf(x)\leq
\tT^n (D\log(x)+q_\vareps(x)+\vareps).
\ee
Corollary \ref{fstarcor} shows that $\{\tT^n(D\log(x)+q_\vareps(x)\,)\}$
converges to $\beta_\vareps f^*$ uniformly on compact subsets
of $(0,2]$, where $\beta_\vareps>0$ is a constant, and 
$\lim_{n\to\infty} \tT 1=w\, f^*$ for $w=-\w_1$.
For $0<y<2$ there exists $n_{\vareps,y}$ such that $j\geq n$ implies
\be
  \sup_{y\leq x\leq2} |\tT^j (D\log(x)+q_\vareps(x)\,)-\beta_\vareps
  f^*(x)|<\vareps,
  \quad \sup_{y\leq x\leq2}|(\tT^j 1)(x) - w f^*(x)|<1.
\ee
Then for $j>n_{\vareps,y}$ and $y\leq x\leq 2$
\be\label{betafstar}
 |\tT^j f(x)-\beta_\vareps f^*(x)|\leq
 \vareps\Big(1+|(\tT^j 1)(x) - w f^*(x)|+wf^*(y)\Big)
 \leq \vareps(\,2+wf^*(y)\,),
\ee
since $f^*(y)\geq f^*(x)$ for $y\leq x\leq 2$.
Thus for $j,k>n_{\vareps,y}$ and $y\leq x\leq 2$,
\be
   |\tT^j f(x)-\tT^k f(x)|\leq \vareps(\,2+wf^*(y)\,),
\ee
so the sequence of functions $\{\tT^n f(x)\}$ is uniformly Cauchy on
each compact subset of $(0,2]$, which implies convergence. It follows
from \eqref{betafstar} this limit must be of the form $\beta^* f^*$,
so we write $\lim_n \tT^n (D\log(x)+q(x))=\beta^* f^*$. By Proposition
\ref{prop:mp}, $\beta^*>0$. Convergence of $\{\tT^n f\}$ to
$\beta^* f^*$ then implies convergence of $\{\hT^n f\}$ to $f^*$.
\end{Proof}

\begin{Corollary}
For $0<x\leq 2$ the sequence $\{\hT^n \delta_x\}$  converges to
$f^*\,dx$.
\end{Corollary}
\begin{Proof}
There exists an integer $j\geq 0$ so that $1<2^j x\leq 2$. Let
$y=2^j x-1$ so that $0<y\leq 1$. Then $T^j \delta_x=\delta_{y+1}$ and 
$T^{j+1} \delta_x=\frac 1y I_{[0,2y]}dx$. There exists an integer
$k\geq 0$ so that $1<2^k y\leq 2$.   Let $z=2^ky-1$.  Then
$T^{j+k+1}\delta_x= \frac1{1+z} (1+ 2I_{[0<x\le2z]}(x)\log\frac{2z}x\;)dx$
which is of the
form $D\log(x)+g(x)$ with $D$ constant and $g$ continuous.
\end{Proof}

% \emph{Acknowledgement.}   The authors wish to thank the anonymous referee for
% corrections and suggested improvements to this paper.

\bibliographystyle{acm}
\def\cprime{$'$}

%\bibliography{/home/mackey/research/bib}

\end{document}